\documentclass[12pt]{article}

\usepackage[english]{babel}
\usepackage[latin1]{inputenc}
\usepackage{times,amsmath,amsthm,amssymb,amsfonts}

\newtheorem{thm}{Theorem}
\newtheorem{lem}[thm]{Lemma}
\newtheorem{cor}[thm]{Corollary}
\newtheorem{exm}[thm]{Example}
\newtheorem{prop}[thm]{Proposition}

\begin{document}


\centerline{\LARGE\bf On the~peripheral spectrum of positive elements}

\medskip

\begin{center}\large \fbox{\it\kern1pt Egor A. Alekhno}\\[9pt]
  Belarusian State University, Minsk, Belarus
\end{center}

\bigskip

{\small \noindent\textbf{Abstract.} Let $A$ be an~ordered Banach algebra with a~unit ${\bf e}$ and
a~cone $A^+$. An~element $p$ of $A$ is said to be an~order idempotent if $p^2=p$ and $0\le p\le{\bf
e}$. An~element $a\in A^+$ is said to be irreducible if the~relation $({\bf e}-p)ap=0$, where
$p$~is an~order idempotent, implies $p=0$ or $p={\bf e}$. For an~arbitrary element $a$ of~$A$
the~peripheral spectrum $\sigma_{\rm per}(a)$ of $a$ is the~set $\sigma_{\rm
per}(a)=\{\lambda\in{\Bbb C}:|\lambda|=r(a)\}$, where $r(a)$ is the~spectral radius of $a$. We
investigate properties of the~peripheral spectrum of an~irreducible element $a$. The~conditions
under which $\sigma_{\rm per}(a)$ contains of or coincides with $r(a)H_m$, where $H_m$ is the~group
of all $m^{\rm th}$ roots of unity, and the~spectrum $\sigma(a)$ is invariant under the~rotation on
angle $\frac{2\pi}m$ for some $m\in{\Bbb N}$, are given. The~correlation between these results and
the~existence of a~cyclic form of $a$ is considered. The~conditions under which $a$ is primitive,
i.e., $\sigma_{\rm per}(a)=\{r(a)\}$, are studied. The~necessary assumptions on the~algebra $A$
which imply the~validity of these results, are discussed. In~particular, the~Lotz-Schaefer axiom is
introduced and finite-rank elements of $A$ are defined. Other approaches to the~notions of
irreducibility and primitivity are discussed. The~conditions under which the~inequalities $0\le
b<a$ imply $r(b)<r(a)$, are studied. The closedness of the~center $A_{\bf e}$, i.e., of~the~order
ideal generated by ${\bf e}$ in $A$, is proved.}

\bigskip

\smallskip

\medskip

\section{Introduction and preliminaries}

Let $A$ be a~(complex) Banach algebra with an~algebraic unit ${\bf e}$ and $A^+$ a~(closed, convex)
cone in $A$. As usual, for elements $a, b \in A$ the~symbol $a\ge b$ (or $b\le a$) means $a - b \in
A^+$. Under this ordering, $A$ is an~ordered linear space. From the~definition of the~cone, it
follows that the~inequalities $a\ge b$ and $b\ge a$ imply $a = b$ for all $a, b \in A$ and $\alpha
x + \beta y\ge 0$ for all elements $x, y\in A^+$ and all scalars $\alpha, \beta\in{\Bbb R}^+$.
The~elements of $A^+$ are called {\it positive}. If ${\bf e}\ge 0$ and the~inequalities $a, b \ge
0$ imply $ab\ge 0$ then $A$ is called an~{\it ordered Banach algebra}. An~important example of
an~ordered Banach algebra is the~algebra of (linear, bounded) operators on an~ordered Banach space.
Namely, if $E$ is an~ordered Banach space with (closed) cone $E^+$ then the~algebra $B(E)$ of all
operators on $E$ is an~ordered Banach algebra under the~natural order if and only if the~linear
space $E^+ - E^+$ is dense in $E$. In~particular, if $E$ is a~(complex) Banach lattice then
the~algebra $B(E)$ is an~ordered Banach algebra.

The~study of ordered Banach algebras was initiated in \mbox{\cite{MR, RR}}. In these papers and in
a~number of subsequent ones the~main emphasis was on the~study of spectral properties of positive
elements. Nevertheless, in spite of the~considerable progress in this direction, several aspects of
the~theory have received almost no attention. In~particular, the~properties of the~peripheral
spectrum of positive elements in ordered Banach algebras were not enough studied. The~main purpose
of this note is to take a~step in this direction.
\smallskip

We recall some well-known notions which will be necessary later on. Let $B$ be a~Banach algebra
with a~unit. The~{\it spectrum} \mbox{\cite[p.~19]{BD}} of an~element $b$ of $B$ is the~set
$$\sigma(b;B)=\{\lambda\in{\Bbb C} : \lambda-b \ \text{is an~invertible element of} \ B\}.$$
When no confusion can occur, we write $\sigma(b)$ to denote $\sigma(b;B)$. The~{\it resolvent
set}~$\rho(b;B)$ \mbox{\cite[p.~245]{AbrAl}} of $b$ is the~complement of the~spectrum, i.e.,
$\rho(b;B)={\Bbb C}\setminus\sigma(b;B)$. As~is well known, the~spectrum $\sigma(b;B)$ is
a~non-empty compact set (we assume $B\neq\{0\}$). The~{\it spectral radius} $r(b)$
\mbox{\cite[p.~23]{BD}} of an~element $b$ is defined via the~formula $r(b)=\max{\{|\lambda| :
\lambda\in\sigma(b;B)\}}$. By the~Gelfand formula \mbox{\cite[pp. 11,~23]{BD}}, the~equality
$r(b)=\lim\limits_{n\to\infty}\|b^n\|_B^{\frac1n}$ holds. The~{\it peripheral spectrum}
$\sigma_{\rm per}(b;B)$ (see, e.g., \mbox{\cite[Section~4]{MR}}) of an~element $b$ is the~set
$$\sigma_{\rm per}(b;B)=\{\lambda\in\sigma(b;B) : |\lambda|=r(b)\}.$$
Obviously, $\sigma_{\rm per}(b;B)$ is also a~non-empty compact set. Next, the~{\it resolvent
function} \mbox{\cite[p.~245]{AbrAl}} $R(\cdot,b):\rho(b;B)\to B$ of $b$ is defined by
$R(\lambda,b)=(\lambda-b)^{-1}$ and is an~analytic function on the~open set $\rho(b;B)$.
In~particular, if $\lambda_0$ is an~isolated point of $\sigma(b;B)$ then in a~sufficiently small
punctured neighbourhood of this point a~Laurent series expansion
\begin{equation}\label{1}
R(\lambda,b)=\sum\limits_{j=-\infty}^{+\infty}b_{\lambda_0,j}(\lambda-\lambda_0)^j
\end{equation}
holds, where $b_{\lambda_0,j}\in B$ for all $j\in{\Bbb Z}$. Of course, this expansion also holds
when the~point $\lambda_0\notin\sigma(b;B)$; in this case, $b_{\lambda_0,j}=0$ for $j\le-1$.
If~$\lambda_0=r(a)$, we will simply write $b_j$ instead of $b_{\lambda_0,j}$. A~point $\lambda_0$
is said to be a~{\it pole} \mbox{\cite[p.~264]{AbrAl}} of the~resolvent $R(\cdot,b)$ of
order~$k\in{\Bbb N}$ whenever $\lambda_0$ is an~isolated point of $\sigma(b;B)$,
$b_{\lambda_0,j}=0$ for $j<-k$, and $b_{\lambda_0,-k}\neq0$; in this case, the~identities
$bb_{\lambda_0,-k}=b_{\lambda_0,-k}b=\lambda_0b_{\lambda_0,-k}$ hold. If $k=1$ then $\lambda_0$ is
called a~{\it simple pole}.
\smallskip

The~considerable progress in the~study of the~peripheral spectrum was attained in the~case of
ordered spaces. The~results concerning to the~spectral theory of non-negative matrices, i.e., to
the~spectral theory in finite dimensional (Archimedean) Riesz spaces, and, in~particular, to
properties of the~peripheral spectrum of non-negative matrices can be found in \mbox{\cite[Chapter
1]{Sch}} and \mbox{\cite[Chapter 8]{AbrAl}}. Properties of the~peripheral spectrum of positive
operators on Banach lattices are considered, e.g., in \cite[Chapter~V and, in~particular, Sections
V.4 and~V.5]{Sch} (see also \mbox{\cite{Al,AbrAl,Gl}}). Some results concerning to the~peripheral
spectrum of positive elements in ordered Banach algebras can be found in \mbox{\cite[Section
4]{MR}}. Below, as far as necessary, we recall some results in these directions. Now we mention
only the~following result~\cite{RR}. If $A$ is an~ordered Banach algebra such that the~spectral
radius function in $A$ is monotone, i.e., the~inequalities $0\le a\le b$ in~$A$ imply $r(a)\le
r(b)$, then $r(a)\in\sigma_{\rm per}(a)$; in particular, this is true when the~cone $A^+$ is normal
and, hence, when the~ordered linear space $A$ is Dedekind complete.
\medskip

The~paper is organized as follows. In the~second section properties of the~peripheral spectrum of
irreducible elements are studied (see, in~particular, Theorem~\ref{thm12}) and some assumptions on
the~algebra $A$ which allow obtaining nice spectral properties, are considered. In~the~next
section, using a~special assumption (the~Lotz-Schaefer axiom), the~results of the~preceding section
will be made more precisely and, moreover, the~notion of finite-rank elements in ordered Banach
algebras and the~conditions of the~primitivity of irreducible elements are discussed. Another
approach to the~notion of irreducibility is considered in the~fourth section. In the~next section,
conditions under which the~inequa\-lities $0\le b<a$ imply $r(b)<r(a)$ are studied. In~the~last
section, the~closedness of the~center $A_{\bf e}$ of an~ordered Banach algebra $A$ is proved.
\medskip

For any unexplained terminology, notations, and elementary properties of ordered Banach spaces, we
refer the~reader to~\cite{AlT}. For information on the~theory of Riesz spaces, Banach lattices, and
operators on these spaces, we suggest~\mbox{\cite{AbrAl, AlB}} (see also \cite{Sch}). More details
on elementary properties of Banach algebras can be found in~\cite{BD} (see also~\cite{BMSmW}).
Throughout the~note, unless stated otherwise, $A$ will stand for an~arbitrary ordered Banach
algebra with a~unit ${\bf e}\neq0$.

\section{The~peripheral spectrum of irreducible elements}

An~important result in the~spectral theory of positive elements is the~{\it theorem about
the~Frobenius normal form}. For the~case of matrices it means that (see, e.g., \mbox{\cite[p.
31]{Sch}}) an~arbitrary non-negative matrix $A$ can be transformed into a~block lower-triangular
matrix
$$\left(\begin{array}{cccc} A_{11} & 0 & \ldots & 0 \\ A_{21} & A_{22} & \ldots & 0 \\ \ldots &
\ldots & \ldots & \ldots \\ A_{n1} & A_{n2} & \ldots & A_{nn}
\end{array}\right)$$
using a simultaneous permutation of rows and columns, where the~matrices $A_{11},\ldots,A_{nn}$ are
irreducible. We recall that a~matrix $D$ is called {\it irreducible} \mbox{\cite[p. 19]{Sch}}
whenever it can not be transformed to the~form $D=\left(\begin{smallmatrix} D_{11} & 0 \\ D_{21} &
D_{22}\end{smallmatrix}\right)$ using a~simultaneous permutation of rows and columns, where
$D_{ii}\neq D$ for $i=1,2$ (the~zero $1\times1$ matrix is irreducible). Obviously,
$\sigma(A)=\bigcup\limits_{i=1}^n\sigma(A_{ii})$. Therefore, we obtain the~dependence of
the~spectrum of an~arbitrary non-negative matrix $A$ on the~spectrum of irreducible matrices.

Now we consider the~case of an~ordered Banach algebra $A$. An~element $p \in A$ is called an~{\it
order idempotent} \cite{Al2} if $0\le p\le{\bf e}$ and $p^2 = p$. Under the~partial ordering
induced by~$A$, the~set of all order idempotents ${\bf OI}(A)$ of $A$ is a~Boolean algebra and its
lattice operations satisfy the~identities $p\wedge q = pq$ and $p\vee q = p + q - pq$ for all $p, q
\in {\bf OI}(A)$ (see~\cite{Al2}). For $p \in {\bf OI}(A)$ and $a\in A$, we put $p^{\rm d}={\bf
e}-p$ and $a_p=pap$. Obviously, $p^{\rm d}\in{\bf OI}(A)$. A~positive element $a\in A$ is said to
be {\it order continuous} \cite{Al2} if $p_\alpha a\downarrow0$ and $ap_\alpha\downarrow0$ in~$A$
whenever $p_\alpha\downarrow0$ in ${\bf OI}(A)$.\footnote{In~\cite{Al3} the~notion of order
continuity was extended to the~case of an~arbitrary element $a\in A$, not necessarily positive.
However, for our purposes, we can limit oneself to this special situation only.} The~collection of
all order continuous element of $A$ will be denoted by $A_{\rm n}$. The~algebra $A$ is called {\it
order regular} if $A_{\rm n}$ is a~subsemi-group of $A$, i.e., $a,b\in A_{\rm n}$ for every $a,b\in
A_{\rm n}$. Next, an~order idempotent $p\in A$ is called $a$-{\it invariant}~\cite{Al2}, where
$a\in A$, if $p^{\rm d}ap=0$. A~positive element $a\in A$ is said to be {\it
irreducible}~\cite{Al2} whenever $a$ has no non-trivial (i.e., $p\neq0$ and $p\neq{\bf e}$)
invariant order idempotents; all other elements of $A$ are called {\it reducible}. An~element $a\in
A^+$ is said to be {\it irreducible with respect to an~order idempotent} $p\in{\bf OI}(A)$ whenever
there exists no $q\in{\bf OI}(A)$ such that $0<q<p$ and $(p-q)aq=0$. An~element $b\in A$ is called
a~{\it block}~\cite{Al2} of an~element $a\in A^+$ if there exists \mbox{$a$-invariant} order
idempotents $p_1$~and~$p_2$ satisfying the~relations $p_2<p_1$ and $b=a_{p_2p_1^{\rm d}}$.
An~element $b$ is called a~{\it spectral block} of an~element $a\in A^+$~\cite{Al2} if $b$ is
a~block of $a$ and $r(b)=r(a)$. An~element $a\in A^+$ is said to {\it have the~Frobenius normal
form}~\cite{Al2} if there exists \mbox{$a$-invariant} order idempotents $p_0,p_1,\ldots,p_n$, which
determine this form, such that ${\bf e}=p_n\ge p_{n-1}\ge\ldots\ge p_0=0$ and if the~relation
$r(a_{p_ip_{i-1}^{\rm d}})=r(a)$ holds for some $i=1,\ldots,n$ then $a_{p_ip_{i-1}^{\rm d}}$ is
irreducible with respect to $p_ip_{i-1}^{\rm d}$. In this case, we have the~inclusion~\cite{Al2}
\begin{equation}\label{2}
\sigma(a)\subseteq\bigcup\limits_{i=1}^n\sigma(a_{p_ip_{i-1}^{\rm d}}).
\end{equation}
An~order continuous element $a\in A$ is said to be {\it spectrally order continuous}~\cite{Al2} if
for every spectral block $b$ of $a$ the~condition that $r(a)$ is a~pole of $R(\cdot,b)$ of order
$k$ implies the~order continuity of the~element $b_{-k}$ (see~(\ref{1})). The~spectral radius
$r(a)$ is called a~{\it finite-rank pole} (abbreviated an~$f$-{\it pole}) of the~resolvent
$R(\cdot,a)$ of a~positive element $a\in A$ if the~inequalities $0\le b\le a$ imply $r(b)\le r(a)$
and if $r(b)=r(a)$ then $r(a)$ is a~pole of $R(\cdot,b)$.

The~theorem about the~Frobenius normal form has the~following form~\cite{Al3}:

\begin{thm}\label{thm1}\
Let an~ordered Banach algebra $A$ be Dedekind complete and let an~ele\-ment~$a$ of $A$ be
spectrally order continuous. Assume that $r(a)>0$ and $r(a)$ is an~$f$-pole of $R(\cdot,a)$. Then
the~element $a$ has the~Frobenius normal form.
\end{thm}

\begin{exm}\label{exm2}\
{\rm {\bf (a)} Let us consider the~ordered Banach algebra $B(E)$, where $E$ is an~arbitrary Banach
lattice. If $T\in B(E)$ then, on the~one hand, the~notion of order continuity can be defined for
the~{\it operator} $T$ on $E$, i.e., \mbox{\cite[p. 46]{AlB}} if a~net
$x_\alpha\stackrel{o}{\longrightarrow}0$ then $Tx_\alpha\stackrel{o}{\longrightarrow}0$.
The~collection of all order continuous operators on $E$ is denoted by $B_{\rm n}(E)$. On the~other
hand, the~notion of order continuity can be defined for $T$ as of an~{\it element} in the~ordered
Banach algebra $B(E)$ which is defined only for a~positive operator. These two notions may differ.
Moreover~\cite{Al3}, the~inclusion $(B_{\rm n}(E))^+\subseteq(B(E))_{\rm n}$ and the~converse one
do not hold in general. Nevertheless, in~the~case of a~Dedekind complete Banach lattice $E$ these
two notions coincide. i.e., we have $(B_{\rm n}(E))^+=(B(E))_{\rm n}$, and, in~particular,
the~algebra $B(E)$ is order regular.

A~positive operator $T$ on a~Banach lattice $E$ is called ({\it band}) {\it
irreducible}~\mbox{\cite[p. 349]{AbrAl}} if $T$~has no non-trivial invariant bands. Obviously,
${\bf OI}(B(E))$ is the~collection of all order projections on $E$. If $E$ is Dedekind complete
then every band $B$ in $E$ is a~projection band and we have the~one-to-one correspondence between
the~set of all bands in $E$ and the~set of all order projections on $E$. Therefore, an~operator $T$
on a~Dedekind complete Banach lattice~$E$ is an~irreducible operator if and only if $T$ is
an~irreducible element of the~ordered Banach algebra $B(E)$. Next, if $E$ is a~Dedekind complete
Banach lattice, the~Lorenz seminorm $\|x\|_L=\inf{\{\sup{\|y_\alpha\|_E:0\le
y_\alpha\uparrow|x|}\}}$ is a~norm on $E$ (e.g., the~order continuous dual $E_{\rm n}^\sim$
separates points of $E$ or $E$ is an~$AM$-space with an~order unit), and $T$ is a~positive order
continuous operator on $E$ with $r(T)>0$ then \cite{Al2} $r(T)$~is an~$f$-pole of the~resolvent
$R(\cdot,T)$ if and only if $r(T)$ is a~finite-rank pole of $R(\cdot,T)$, i.e., $r(T)$ is a~pole of
$R(\cdot,T)$ and the~residue $T_{-1}$ of $R(\cdot,T)$ at $r(T)$ is a~finite-rank operator. Various
assumptions under which an~operator $T$ has the~Frobenius normal form in a~special case of
the~algebra $B(E)$ can be found in~\cite{Al2} .

{\bf (b)} If $E$ is an~ordered linear space and an~element $x\in E^+$ then the~{\it order ideal
$E_x$ generated by} $x$ is the~set \mbox{\cite[p. 103]{AlT}}
\begin{equation}\label{9}
E_x=\{y\in E:-\lambda x\le y\le\lambda x \ \text{for some} \ \lambda\in{\Bbb R}^+\}.
\end{equation}
Under the~algebraic operations and the~ordering induced be $E$, $E_x$ is a~real ordered linear
space satisfying $E_x\subseteq E^+-E^+$. If $A$ is an~ordered Banach algebra and $b\in A_{\rm n}$
then $A_b^+\subseteq A_{\rm n}$. The~order ideal $A_{\bf e}$ is called \cite{Al3} the~{\it center}
of $A$. As will be shown in Section~\ref{sec6}, $A_{\bf e}$ is closed in~$A$ and, hence, is a~{\it
real} ordered Banach algebra. Again, if ${\bf e}\in A_{\rm n}$ then $A_{\bf e}^+\subseteq A_{\rm
n}$. However, in general, the~inclusion ${\bf e}\in A_{\rm n}$ does not hold (see~\cite{Al3}).
Nevertheless, in every case, the~algebra $A_{\bf e}$ is order regular. Indeed, let $a\in(A_{\bf
e})_{\rm n}$ and $b\in A_{\bf e}^+$. If $p_\alpha\downarrow0$ in ${\bf OI}(A_{\bf e})$ then for
some $\lambda\ge0$, we have $0\le p_\alpha ab\le\lambda p_\alpha a\downarrow0$ in $A_{\bf e}$;
analogously, $abp_\alpha\downarrow0$. Thus, $ab,ba\in(A_{\bf e})_{\rm n}$.

The~author does not know an~example of an~ordered Banach algebra which is not order
regular.}\hfill$\Box$
\end{exm}

It follows from the~theorem about the~Frobenius normal form (see Theorem~\ref{thm1}) and
the~inclusion~(\ref{2}) that the~spectrum of a~wide class of positive elements of an~ordered Banach
algebra~$A$ is determined by the~spectra of irreducible elements. The~next result and
Corollary~\ref{cor5} make more precisely this result for the~case of the~peripheral spectrum.

\begin{lem}\label{lem3}\
Let $B$ be a~Banach algebra with a~unit ${\bf u}$. Let elements $b,p_0,p_1,\ldots,p_n\in B$ with
$n\in{\Bbb N}$ satisfy $({\bf u}-p_j)bp_j=0$, $p_0=0$, $p_n={\bf u}$, and
$p_ip_j=p_{\min{\{i,j\}}}$ for all $i,j=0,\ldots,n$. Then $r(b_{q_j})\le r(b)$ and
$$\sigma_{\rm per}(b)=\bigcup\{\sigma_{\rm per}(b_{q_j}):r(b_{q_j})=r(b)\},$$
where $q_j=p_j-p_{j-1}$ and $b_{q_j}=q_jbq_j$ for $j=1,\ldots,n$.
\end{lem}

{\bf Proof.} We mention first that for an~arbitrary scalar $\lambda$ belonging to the~unbounded
connected component $\rho_\infty(b)$ of the~resolvent set $\rho(b)$ the~identity $({\bf
u}-p_j)R(\lambda,b)p_j=0$ holds for all $j=0,\ldots,n$. Indeed, using an~elementary induction, we
get $({\bf u}-p_j)b^np_j=0$ for all $n\in{\Bbb N}$. Whence, in view of the~expansion
$R(\lambda,b)=\frac1\lambda+\frac1{\lambda^2}b+\ldots$ with $\lambda>r(b)$, we have the~identity
$({\bf u}-p_j)R(\lambda,b)p_j=0$. Taking the~unique extension of the~analytic function $({\bf
u}-p_j)R(\lambda,b)p_j$ to $\lambda\in\rho_\infty(b)$, we conclude that the~last identity holds for
all $\lambda\in\rho_\infty(b)$; in~particular,
\begin{equation}\label{3}
({\bf u}-p_j)R(\lambda,b)q_j=0.
\end{equation}

Now let us consider a~non-zero scalar $\lambda\in\rho_\infty(b)$. The~proof will be completed if we
will check the~invertibility of the~element $\lambda-b_{q_j}$ for every $j=1,\ldots,n$. To this
end, we define the~element $z$ via the~formula $z=\frac1\lambda({\bf
u}-p_j)+q_jR(\lambda,b)q_j+\frac1\lambda p_{j-1}$. Using (\ref{3}) and the~equality
$q_jbp_{j-1}=0$, we have
$$(\lambda-b_{q_j})z={\bf u}-p_j+(\lambda-b_{q_j})q_jR(\lambda,b)q_j+p_{j-1}=
{\bf u}-q_j+q_j(\lambda-b_{q_j})R(\lambda,b)q_j=$$
$$={\bf u}-q_j+q_j(\lambda-b+b({\bf u}-q_j))R(\lambda,b)q_j=
{\bf u}+q_jb({\bf u}-q_j)R(\lambda,b)q_j=$$
$$={\bf u}+q_jb({\bf u}-p_j+p_j-q_j)R(\lambda,b)q_j=
{\bf u}+(q_jb({\bf u}-p_j)+q_jbp_{j-1})R(\lambda,b)q_j=$$
$$={\bf u}+q_jb({\bf u}-p_j)R(\lambda,b)q_j={\bf u}.$$
Analogously, $z(\lambda-b_{q_j})={\bf u}$.\hfill$\Box$
\smallskip

As the~next example shows, the~inclusion~(\ref{2}) cannot make more precisely in general.

\begin{exm}\
{\rm If $a\in A$, $q$ is an~idempotent of $A$, and $q\neq {\bf e}$ then $0\in\sigma(a_q)$. Thus,
the~inclusion~(\ref{2}) is proper in general. We shall show that the~identity
$\sigma(a)\setminus\{0\}=(\sigma(a_{p^{\rm d}})\cup\sigma(a_p))\setminus\{0\}$ also does not hold,
where $p$ is an~\mbox{$a$-invariant} order idempotent. To see this, we consider the~space
$\ell_\infty$ of all bounded sequences and define the~operator $S$ on $\ell_\infty$ via the~formula
$$Sx=(x_1,0,0,0,x_3,0,0,0,x_5,\ldots)+(0,x_2,x_4,x_6,0,x_8,x_{10},x_{12},0,\ldots),$$
where $x=(x_1,x_2,\ldots)\in\ell_\infty$. As is easy to see, $S$ is invertible. Put $T=I+S$.
Obviously, $1\notin\sigma(T)$. The~band $B=\{x\in\ell_\infty:x_{2k}=0 \ {\rm for \ all} \ k\in{\Bbb
N}\}$ is \mbox{$T$-invariant} and
$$P_BTP_Bx=(x_1,0,x_3,0,x_5,0,x_7,0,x_9,\ldots)+(x_1,0,0,0,x_3,0,0,0,x_5,\ldots),$$
where $P_B$ is the~order projection on $B$. In~particular, $1\in\sigma(P_BTP_B)$.}\hfill$\Box$
\end{exm}

\begin{cor}\label{cor5}\
Let order idempotents $p_0,p_1,\ldots,p_n$ determine the~Frobenius normal form of an~element $a\in
A$. Then
$$\sigma_{\rm per}(a)=\bigcup\{\sigma_{\rm per}(a_{q_j}):r(a_{q_j})=r(a) \ \text{and} \ a_{q_j} \
\text{is irreducible with respect to} \ q_j\},$$ where $q_j=p_jp_{j-1}^{\rm d}$ for $j=1,\ldots,n$.
\end{cor}

Keeping the~preceding corollary in mind, we now turn to the~study of the~peripheral spectrum of
irreducible elements. Firstly, let us recall some spectral properties of irreducible elements. It
should be mentioned at once that there exists an~ordered Banach algebra $A$ such that every
positive element of $A$ is irreducible; the~latter is equivalent to the~identity ${\bf
OI}(A)=\{0,{\bf e}\}$. For example, the~ordered Banach algebra $A=A_0\otimes{\Bbb C}$ obtained from
an~ordered Banach algebra~$A_0$ by adjoining a~unit or the~algebras $C(K)$ of  all continuous
functions on $K$ and $B(C(K))$, where $K$ is a~connected (Hausdorff) compact. Therefore, in
general, one cannot expect any distinctive spectral properties of irreducible elements and, hence,
a~special class of ordered Banach algebras should be distinguished. We shall say that an~ordered
Banach algebra $A$ has a~{\it disjunctive product}~\cite{Al2} if for any $a,b\in A_{\rm n}$ with
$ab=0$ there exists an~element $p\in{\bf OI}(A)$ satisfying $ap=p^{\rm d}b=0$. The~algebra $B(E)$,
where $E$ is a~Dedekind complete Banach lattice, has~\cite{Al2} a~disjunctive product. Next,
an~element $a\in A$ is said to be {\it algebraically strictly positive}, in symbols $a\gg0$,
whenever $p_1ap_2>0$ for all $0<p_1,p_2\in{\bf OI}(A)$. The~next result holds~\cite{Al2}.

\begin{thm}\label{thm6}\
Let $A$ be an~ordered algebra with a~disjunctive product and with the~Boo\-lean algebra ${\bf
OI}(A)$ Dedekind complete. Assume that a~non-zero element $a\in A$ is order continuous, that $r(a)$
is a~pole of $R(\cdot,a)$ of order $k$, and that $a_{-k}$ is also order continuous. If the~element
$a$ is irreducible then the~following statements hold:
\begin{description}
\item[(a)] The~spectral radius $r(a)>0$;
\item[(b)] The~point $r(a)$ is a~simple pole of $R(\cdot,a)$;
\item[(c)] The~residue $a_{-1}\gg0$ and the~resolvent $R(\lambda,a)\gg0$ for all
$\lambda>r(a)$;
\item[(d)] If $0\le b<a$ and if some element $c\in A_{\rm n}$ satisfies $0< r(b)c\le bc$ then $r(b)<r(a)$.
\end{description}
\end{thm}

The~following lemma and Corollary~\ref{cor8} show that, under additional assumptions on~$A$, the
condition about the~order continuity of the~coefficient $a_{-k}$ of the~Laurent series expansion of
$R(\,\cdot\,, a)$ around $r(a)$ in the~preceding theorem can be rejected.

\begin{lem}\
Let an~ordered Banach algebra $A$ be order regular and let $a\in A_{\rm n}$. Then the~resolvent
$R(\lambda,a)\in A_{\rm n}$ for all $\lambda>r(a)$.
\end{lem}

{\bf Proof.} For an~arbitrary number $n\in{\Bbb N}$, we define the~elements $b_n$ and $c_n$ as
follows
$$
 b_n=\frac1\lambda{\bf e}+\ldots+\frac1{\lambda^n}a^{n-1} \ \ \text{and} \ \
 c_n=\frac1{\lambda^{n+1}}a^n+\frac1{\lambda^{n+2}}a^{n+1}+\ldots.
$$
Obviously, $b_n+c_n=R(\lambda,a)$, $c_n\to0$ as $n\to\infty$, and, since $A$ is order regular,
$b_n\in A_{\rm n}$. Let $p_\alpha\downarrow0$ in ${\bf OI}(A)$ and let $p_\alpha R(\lambda,a)\ge x$
with $x\in A$. Fix an~index $\alpha_0$. For every $\alpha\ge\alpha_0$, the~inequality $x\le
p_\alpha b_n+p_{\alpha_0}c_n$ holds and, hence, $x-p_{\alpha_0}c_n\le p_\alpha
b_n\downarrow_{\alpha\ge\alpha_0}0$. Therefore, $x\le p_{\alpha_0}c_n\to0$ as $n\to\infty$.
Finally, $x\le0$ or $p_\alpha R(\lambda,a)\downarrow0$. Analogously, $R(\lambda,a)p_\alpha
\downarrow0$.\hfill$\Box$

\begin{cor}\label{cor8}\
Let an~ordered Banach algebra $A$ be order regular and let $A_{\rm n}$ be closed in~$A$. Let
an~element $a\in A_{\rm n}$ such that $r(a)$ is a~pole of the~resolvent $R(\cdot,a)$ of order~$k$.
Then $a_{-k}\in A_{\rm n}$.
\end{cor}

{\bf Proof.} The~identity $a_{-k}=\lim\limits_{\lambda\downarrow r(a)}(\lambda-r(a))^kR(\lambda,a)$
holds. It only remains to recall the~preceding lemma and the~closedness of $A_{\rm n}$.\hfill$\Box$
\smallskip

If $E$ is a~Dedekind complete Banach lattice such that the~Lorenz seminorm~\mbox{$\|\cdot\|_L$} is
a~norm on $E$ (see Example~\ref{exm2}{\bf (a)}) then \cite{KW} the~set $(B(E))_{\rm n}$ is closed
in $B(E)$. Consequently, the~wide class of ordered Banach algebras of the~form $B(E)$
auto\-ma\-ti\-cally satisfies the~assumption about the~closedness of $A_{\rm n}$.

Thus, it follows from Theorem~\ref{thm6} and the~preceding corollary that, under the~next
assumptions on $A$, order continuous irreducible elements have nice spectral properties:
\begin{description}
\item[{\bf (A$_1$)}] An~ordered Banach algebra $A$ is order regular and has a~disjunctive product,
the~Boolean algebra ${\bf OI}(A)$ is Dedekind complete, and the~set $A_{\rm n}$ is closed in $A$.
\end{description}

If $E$ is a~Dedekind complete Banach lattice and $P$ is a~non-zero order continuous projection on
$E$ such that $P\gg0$ as an~element of $B(E)$ then \cite{Al2} $\dim{R(P)}=1$, where $R(P)$ is
the~range of the~operator $P$. In~particular, in view of part~{\bf (c)} of Theorem~\ref{thm6}, ``as
a rule", the~residue $T_{-1}$ of the~resolvent $R(\cdot,T)$ of an~irreducible operator $T$ at
$r(T)$ is a~rank-one operator. On the~other hand, as is well known, if $Z$ is a~Banach space and
a~projection $Q\in B(Z)$ then $Q$ is a~rank-one operator if and only if $Q$ is a~minimal idempotent
of the~algebra $B(Z)$, i.e., $QB(Z)Q={\Bbb C}Q$. In the~case of an~arbitrary ordered Banach algebra
$A$ and of an~idempotent $b\in A$, the~condition $b\gg0$ does not imply the~minimality of
the~idempotent $b$. For this reason, we must axiomatize this property and make the~following
assumption:
\begin{description}
\item[{\bf (A$_2$)}] Every algebraically strictly positive idempotent $b$ of $A$ is minimal.
\end{description}

As usual, through $L_a$ and $R_a$, we will denote the~operators on an~algebra $B$ defined by
\begin{equation}\label{7}
L_ab=ab \ \ \text{and} \ \ R_ab=ba,
\end{equation}
where $a,b\in B$, and through $N(S)$, we will denote the~null space of the~operator $S$ acting
between two linear spaces.

\begin{prop}\label{prop9}\
An~ordered Banach algebra $A$ satisfies Axiom~{\bf (A$_2$)} if and only if for an~arbitrary
idempotent $b\gg0$ of the~following identity holds
\begin{equation}\label{4}
\dim{N(I-L_b)\cap N(I-R_b)}=1.
\end{equation}
\end{prop}

In general, the~next result is true: {\it If $B$ is a~Banach algebra and an~idempotent $b\in B$
then $b$ is minimal if and only if {\rm (\ref{4})} holds}.
\smallskip

{\bf Proof.} Let $A$ satisfy Axiom~{\bf (A$_2$)} and let $a\in A$ such that $a=ba=ab$. Then
$a=bab=\lambda b$ for some $\lambda\in{\Bbb C}$. For the~converse, if $b^2=b\gg0$ then for every
$a\in A$ the~element $bab$ belongs to $N(I-L_b)\cap N(I-R_b)$. On the~other hand, this space
contains $b$ and, hence, $bab=\mu b$ for some $\mu\in{\Bbb C}$.\hfill$\Box$
\smallskip

For the~study of the~peripheral spectrum of a~positive operator on a~Banach lattice,
the~possibility of the~restriction of a~problem to the~case of operators on the~space $C(K)$ of all
continuous functions on a~compact space $K$ which is simpler for the~study, is important. Recall
that if $E$ is a~Riesz space satisfying Axiom~{\bf (OS)}, i.e., \mbox{\cite[p. 54]{Sch}} if
the~inequalities $0\le z_n\le\lambda_n z$, where $z_n,z\in E$ and
$\lambda=(\lambda_1,\lambda_2,\ldots)\in\ell_1$, imply the~existence of
$\sup\limits_n\sum\limits_{j=1}^nz_j$, then the~order ideal $E_x$ generated by $x$ is, under
the~Minkowski norm $\|\cdot\|_x$ defined by
\begin{equation}\label{28}
\|y\|_x=\inf{\{\lambda\in{\Bbb R}^+:-\lambda x\le y\le\lambda x\}}
\end{equation}
with $y\in E_x$, an~$AM$-space with order unit $x$~\mbox{\cite[p. 102]{Sch}}. Therefore, by
the~Kakutani-Bohnenblust-M.-S.Krein theorem~\mbox{\cite[p. 201]{AlB}}, $E_x$ is lattice isometric
onto a~spa\-ce~$C(K)$ and, moreover, under this isomorphism, $x$ is mapped onto the~constant-one
function~${\rm 1}\!\!{\rm 1}_K$. Evidently, every Dedekind complete Riesz space $E$ satisfies
Axiom~{\bf (OS)}. Now let $E$~be a~Dedekind complete Banach lattice. Recall that an operator $T$ on
$E$ is said to be {\it regular}~\mbox{\cite[p. 12]{AlB}} if it can be written as a~difference of
two positive operators. As is well known, every regular operator is bounded. By
the~Riesz-Kantorovich theorem~\mbox{\cite[p. 14]{AlB}}, the~space $L_{\rm r}(E)$ of all regular
operators on $E$ is a~Dedekind complete Riesz space. Hence, if $T\in L_{\rm r}(E)$ then the~order
ideal $(B(E))_T$ is lattice isometric onto a~space $C(K)$. We axiomatize this property and make
the~next assumption:
\begin{description}
\item[{\bf (A$_3$)}] The~order
ideal $A_b$ generated by the~non-zero element $b\in A^+$ is lattice isomorphic onto a~space $C(K)$
and, under this isomorphism, $x$ is mapped onto~${\rm 1}\!\!{\rm 1}_K$.
\end{description}

Every complex Riesz space $E$ is the~complexification of the~real Riesz space $E_{\Bbb R}$
satisfying Axiom~{\bf (OS)} (see \mbox{\cite[Section~II.11]{Sch}}). In this case, the~algebra
$L(E)$ of all operators on $E$ is isomorphic onto the~complexification of the~algebra $L(E_{\Bbb
R})$ and, in~particular, every $T\in L(E)$ has a~unique decomposition $T=T_1+iT_2$, where $T_j$ are
real maps on $E_{\Bbb C}$, i.e., $T_j(E_{\Bbb R})\subseteq E_{\Bbb R}$ for $j=1,2$. Defining
the~set $A_{\rm r}$ of all {\it regular}~\cite{Al3} elements of an~ordered Banach algebra $A$ by
$A_{\rm r}=A^+-A^+$, we axiomatize this property and make the~next assumption:
\begin{description}
\item[{\bf (A$_4$)}] The~equality $a+ib=0$ with $a,b\in A_{\rm r}$ implies $a=b=0$.
\end{description}

We continue our discussion with two auxiliary results.

\begin{lem}\label{lem10}\
Let $B$ be a~Banach algebra with a~unit, let $b\in B$, and let $m,k\in{\Bbb N}$. If
$\lambda_0\in\sigma(b)$ and the~set
$\{\lambda_0,\lambda_0e^{i\frac{2\pi}m},\ldots,\lambda_0e^{i\frac{2\pi}m(m-1)}\}\cap\sigma(b)$
consists entirely of poles of $R(\cdot,b)$ of orders which are not greater than $k$, then
$\lambda_0^m$ is a~pole of $R(\cdot,b^m)$ of order which is not greater than~$km$.
\end{lem}

{\bf Proof.} Put $\omega_j=e^{i\frac{2\pi}mj}$ with $j=0,\ldots,m-1$. We claim first that
$\lambda_0^m$ is an~isolated point of $\sigma(b^m)$. Indeed, if a~sequence $\{\xi_n\}$ in
$\sigma(b^m)$ satisfies $\xi_n\neq\lambda_0^m$ for all $n$ and $\xi_n\to\lambda_0^m$ as
$n\to\infty$ then, taking into account the~identity $\sigma(b^m)=f(\sigma(b))$ with $f(z)=z^m$, we
find a~sequence~$\{\mu_n\}$ in $\sigma(b)$ with the~property $\mu_n^m=\xi_n$. Let $\{\mu_{n_r}\}$
be an~arbitrary convergent subsequence of~$\{\mu_n\}$. If $\mu_{n_r}\to\mu_0\in\sigma(b)$ then
$\mu_0^m=\lambda_0^m$ and, hence, $\mu_0=\lambda_0\omega_j$ for some $j=0,\ldots,m-1$. Therefore,
$\mu_0$ is an~isolated point of~$\sigma(b)$. Thus, $\mu_{n_r}=\mu_0$ for sufficiently large $r$ and
for such~$r$, we have $\xi_{n_r}=\lambda_0^m$, a~contradiction.

For arbitrary numbers $\lambda,z\in{\Bbb C}$ the~identity
$\lambda^m-z^m=\prod\limits_{j=0}^{m-1}(\lambda-\omega_jz)$ holds. The~latter implies
\begin{equation}\label{5}
\lambda^m-b^m=\prod\limits_{j=0}^{m-1}(\lambda-\omega_jb).
\end{equation}
If ${\cal U}_{\lambda_0^m}$ is a~punctured neighbourhood of the~point $\lambda_0^m$ satisfying
${\cal U}_{\lambda_0^m}\subseteq\rho(b^m)$ then the~set ${\cal V}=\{\lambda\in{\Bbb
C}:\lambda^m\in{\cal U}_{\lambda_0^m}\}$ is open and the~inclusion ${\cal
V}\subseteq\rho(\omega_jb)$ holds for all $j=0,\ldots,m-1$. Now, using~(\ref{5}), we obtain
\begin{equation}\label{6}
R(\lambda^m,b^m)=\prod\limits_{j=0}^{m-1}R(\lambda,\omega_jb)
\end{equation}
for all $\lambda\in{\cal V}$. Let us consider an~arbitrary sequence $\{\lambda_n\}$ in ${\Bbb C}$
satisfying $\lambda_n\neq\lambda_0^m$ for all $n$ and $\lambda_n\to\lambda_0^m$. It is not
difficult to show that there exists a~sequence $\{\theta_n\}$ in ${\Bbb C}$ with the~properties
$\theta_n\neq\lambda_0^m$, $\theta_n^m=\lambda_n$ for all $n$, and $\theta_n\to\lambda_0$ as
$n\to\infty$. Obviously, $\theta_n\in{\cal V}$ for all sufficiently large $n$. In view of
the~identity~(\ref{6}), we have
$$(\lambda_n-\lambda_0^m)^{km}R(\lambda_n,b^m)=
(\theta_n^m-\lambda_0^m)^{km}R(\theta_n^m,b^m)=
(\theta_n^m-\lambda_0^m)^{km}\prod\limits_{j=0}^{m-1}R(\theta_n,\omega_jb)=$$
$$=\prod\limits_{j=0}^{m-1}(\theta_n^m-\lambda_0^m)^kR(\theta_n,\omega_jb)=
(-1)^{m+1}\prod\limits_{j=0}^{m-1}(\theta_n^m-\lambda_0^m)^kR(\theta_n\omega_j^{-1},b)=$$
$$=(-1)^{m+1}\prod\limits_{j=0}^{m-1}(\theta_n^m-\lambda_0^m)^kR(\theta_n\omega_j,b)=$$
$$=(-1)^{m+1}
\prod\limits_{j=0}^{m-1}\Big(\sum\limits_{l=0}^{m-1}\theta_n^{m-1-l}\lambda_0^l\Big)^k\cdot
\prod\limits_{j=0}^{m-1}(\theta_n-\lambda_0)^kR(\theta_n\omega_j,b)=$$
$$=(-1)^{(m+1)(k+1)}
\prod\limits_{j=0}^{m-1}\Big(\sum\limits_{l=0}^{m-1}\theta_n^{m-1-l}\lambda_0^l\Big)^k
\cdot\prod\limits_{j=0}^{m-1}(\theta_n\omega_j-\lambda_0\omega_j)^kR(\theta_n\omega_j,b)\to 0$$ as
$n\to\infty$.\hfill$\Box$
\medskip

Let $m\in{\Bbb N}$. Consider a~nonempty subset $J$ of the~set $\{1,\ldots,m\}$. Let
$J=\{j_1,\ldots,j_r\}$ with $r=1,\ldots,m$. We define the~shift $J-1$ of $J$ via the~formula
$J-1=\{j_1-1,\ldots,j_r-1\}$, where in the~case of $j_k=1$ for some $k=1,\ldots,r$, we put
$j_k-1=m$. Now, using the~elementary induction, the~set $J-l$ can be defined easily for every
$l\in{\Bbb N}$. Obviously, $J-m=J$.

\begin{lem}\label{lem11}\
Let $m\in{\Bbb N}$ and let $J$ be a~nonempty subset of $\{1,\ldots,m\}$. If $l$ is a~minimal
natural number satisfying $J-l=J$ then $l$ is a~divisor of $m$.
\end{lem}

{\bf Proof.} Evidently, $J-kl=J$ for all $k\in{\Bbb N}$ and $l\le m$. The~representation $m=nl+r$
holds with $n,r\in{\Bbb N}$ and $0\le r<l$. Therefore, $J=J-(nl+r)=J-r$. Whence, using
the~minimality of $l$, we infer $r=0$.\hfill$\Box$
\smallskip

Let be a~Banach algebra with a~unit. As is well known, if an~element $a\in B$ then the~identities
$\sigma(L_a;B(A))=\sigma(R_a;B(A))=\sigma(a;A)$ hold (see~(\ref{7})). We shall say that a~point
$\lambda\in\sigma(a)$ is the~{\it joint eigenvalue} of $a$ whenever there exists a~non-zero element
$c\in B$ satisfying
\begin{equation}\label{8}
ac=ca=\lambda c.
\end{equation}
Now let $A$ be an~ordered Banach algebra, let $a\in A$, and let $b\in A^+$. The~{\it joint
spectrum} of an~element $a$ with respect of $b$ is the set $\sigma_{\rm j}(a;b;A)$ of all complex
numbers $\lambda$ such that there exists an~element $c\in A$ which is not nilpotent,
satisfies~(\ref{8}), and has the~representation in the~form $c=c'+ic''$, where $c'$ and $c''$
belong to the~order ideal $A_b$ (see~(\ref{9})). Again, when no confusion can occur, we write
$\sigma_{\rm j}(a;b)$ to denote $\sigma_{\rm j}(a;b;A)$.
\smallskip

Below, through $H_m$ with $m\in{\Bbb N}$, we will denote the~group of all $m^{\rm th}$ roots of
unity.
\smallskip

Now we are ready to state and to prove the~main result of this paper.

\begin{thm}\label{thm12}\
Let an~ordered Banach algebra $A$ satisfy Axioms {\bf (A$_1$)}-{\bf (A$_4$)}. Let $a$ be
a~non-zero, order continuous, and irreducible element of $A$ such that $r(a)$ is a~pole
of~$R(\cdot,a)$. Let $m\in{\Bbb N}$, $m>1$. The~following statements are equivalent:
\begin{description}
\item[(a)] $a^m$ is reducible;
\item[(b)] For some divisor $m'$ of $m$, $m'>1$, there exist non-zero order idempotents
$p_1,\ldots,p_{m'}$ of $A$ satisfying $\sum\limits_{j=1}^mp_j={\bf e}$, $p_{j'}p_{j''}=0$ for
$j'\neq j''$, and $p_ja=ap_{j+1}$ for all $j=1,\ldots,m'$, where for $j=m'$, we put $j+1=1$;
\item[(c)] For some divisor $m''$ of $m$, $m''>1$, and for all $j=0,\ldots,m''-1$ the~points
$r(a)e^{i\frac{2\pi}mj}$ are poles of $R(\cdot,a)$.
\end{description}

If the~condition~{\bf (b)} holds then the~spectrum $\sigma(a)$ is invariant under the~rotation on
angle~$\frac{2\pi}{m'}$, i.e., $\sigma(a)=\frac{2\pi}{m'}\sigma(a)$.

Moreover, $\sigma_{\rm per}(a)\cap\sigma_{\rm j}(a;a_{-1})=r(a)H_{m_0}$ for some $m_0\in{\Bbb N}$.
\end{thm}

Below, we shall say that an~arbitrary element $a\in A$ has the~{\it cyclic form} whenever there
exist order idempotents $p_1,\ldots,p_{m'}$ with $m'>1$ satisfying the~conditions of part~{\bf
(b)}. In this case, $p_1,\ldots,p_{m'}$ is said to determine the~cyclic form of $a$.
\medskip

{\bf Proof.} {\bf (a)}~$\Longrightarrow$~{\bf (b)} Taking into account the~reducibility of $a^m$,
we find a~non-trivial order idempotent $q_1$ of $A$ satisfying $q_1^{\rm d}a^mq_1=0$. Obviously,
$q_1^{\rm d}aa^{m-1}q_1=0$. Since the~algebra $A$ has a~disjunctive product (Axiom~{\bf (A$_1$)}),
there exists an~order idempotent~$q_2$ of $A$ such that $q_1^{\rm d}aq_2=q_2^{\rm d}a^{m-1}q_1=0$.
Using the~elementary induction, we find $q_1,\ldots,q_m\in{\bf OI}(A)$ satisfying $q_1^{\rm
d}aq_2=q_2^{\rm d}aq_3=\ldots=q_m^{\rm d}aq_1=0$. The~identity
\begin{equation}\label{10}
\prod\limits_{j=1}^mq_j=0
\end{equation}
holds. To see this, keeping the~relation
\begin{equation}\label{11}
aq_j=q_{j-1}aq_j
\end{equation}
for $j=1,\ldots,m$ (for $j=1$, we put $j-1=m$) in mind, we have
$$\Big({\bf e}-\prod\limits_{j=1}^mq_j\Big)a\prod\limits_{j=1}^mq_j=
\Big(q_m-\prod\limits_{j=1}^mq_j\Big)a\prod\limits_{j=1}^mq_j=$$
$$=\Big(q_mq_1-\prod\limits_{j=1}^mq_j\Big)a\prod\limits_{j=1}^mq_j=\ldots=
\Big(q_mq_1\ldots q_{m-1}-\prod\limits_{j=1}^mq_j\Big)a\prod\limits_{j=1}^mq_j=0.$$
The~irreducibility of $a$ and the~relation $q_1\neq{\bf e}$ imply~(\ref{10}).

For $r\in\{1,\ldots,m\}$, we put ${\cal P}_r=\{J\subseteq\{1,\ldots,m\}:{\rm card}\! \ J=r\}$. Now
we assume the~validity of the~identity
\begin{equation}\label{12}
\prod\limits_{k\in J}q_k=0
\end{equation}
for all $J\in{\cal P}_r$, where $r=2,\ldots,m$ and for $r=m$, we obtain~(\ref{10}). Let us show
that in this case either there exists the~required collection of idempotents or
the~identity~(\ref{12}) is valid for all $J\in{\cal P}_{r-1}$. To this end, let (\ref{12}) hold for
all $J\in{\cal P}_r$, $r\ge2$. Then for all $J_1,J_2\in{\cal P}_{r-1}$, $J_1\neq J_2$, we have
\begin{equation}\label{13}
\prod\limits_{j\in J_1}q_j\cdot\prod\limits_{k\in J_2}q_k=0.
\end{equation}
We shall say that two subsets $J_1,J_2\in{\cal P}_{r-1}$ is {\it equivalent} whenever $J_1-l=J_2$
for some $l\in{\Bbb N}$ (see the~remarks before Lemma~\ref{lem11}). As is easy to see, the~relation
introduced above is an~equivalence relation on ${\cal P}_{r-1}$ really. Thus, the~set ${\cal
P}_{r-1}$ is the~union of (disjoint) equivalence classes ${\cal I}_1,\ldots,{\cal I}_t$ with
$t\in{\Bbb N}$. Clearly, if $J\in{\cal I}_s$, where $s\in1,\ldots,t$, then $J-1\in{\cal I}_s$.
By~Lemma~\ref{lem11}, $m_s={\rm card}\! \ {\cal I}_s$ is a~divisor of $m$. If $J\in{\cal I}_s$ then
\begin{equation}\label{14}
{\cal I}_s=\{J,J-1,\ldots,J-(m_s-1)\}
\end{equation}
Fix $s$. In view of (\ref{11}), (\ref{12}), (\ref{13}), and~(\ref{14}), for an~arbitrary set
$J_0\in{\cal I}_s$, we have
$$\Big({\bf e}-\sum\limits_{J\in{\cal I}_s}\prod\limits_{j\in J}q_j\Big)
a\prod\limits_{j\in J_0}q_j=$$
$$=\Big({\bf e}-\sum\limits_{J\in{\cal I}_s}\prod\limits_{j\in J}q_j\Big)
\prod\limits_{k\in J_0-1}q_j\cdot a\prod\limits_{j\in J_0}q_j= \Big(\prod\limits_{k\in
J_0-1}q_j-\prod\limits_{k\in J_0-1}q_j\Big)a\prod\limits_{j\in J_0}q_j=0.$$ Since $J_0$ is
arbitrary, we obtain $\Big({\bf e}-\sum\limits_{J\in{\cal I}_s}\prod\limits_{j\in J}q_j\Big)a
\sum\limits_{J\in{\cal I}_s}\prod\limits_{j\in J_0}q_j=0$. Thus, for every $s\in1,\ldots,t$ either
$\sum\limits_{J\in{\cal I}_s}\prod\limits_{j\in J_0}q_j=0$ or $\sum\limits_{J\in{\cal
I}_s}\prod\limits_{j\in J_0}q_j={\bf e}$. Moreover, either the~former of these equalities hold for
all indexes $s$ or, in view of~(\ref{13}), the~second one holds for the~unique index~$s_0$ (if
$r=m$ then $s_0=t=1$ and $m_{s_0}=m$). We consider the~second case. Let ${\cal
I}_{s_0}=\{J_1,\ldots,J_{m_{s_0}}\}$ and let $J_{j-1}=J_j-1$ for all $j=1,\ldots,m_{s_0}$. Now we
define required order idempotents $p_1,\ldots,p_{m_{s_0}}$ via the~formula $p_j=\prod\limits_{r\in
J_j}q_r$ for $j=1,\ldots,m_{s_0}$. Obviously, $\sum\limits_{j=1}^{m_{s_0}}p_j={\bf e}$ and
$p_{j'}p_{j''}=0$ for $j'\neq j''$. Next, for an~arbitrary index $j=1,\ldots,m_{s_0}$, we have
$$ap_{j+1}=a\prod\limits_{r\in J_j}q_r=
\sum\limits_{J\in{\cal I}_{s_0}}\prod\limits_{r\in J}q_r\cdot a\cdot\prod\limits_{r\in J_{j+1}}q_r
=$$
$$=\Big(\sum\limits_{J\in{\cal I}_{s_0}}\prod\limits_{r\in J}q_r\Big)
\prod\limits_{r\in J_{j+1}-1}q_r\cdot a= \Big(\sum\limits_{J\in{\cal I}_{s_0}}\prod\limits_{r\in
J}q_r\Big)\prod\limits_{k\in J_j}q_k\cdot a= \prod\limits_{r\in J_j}q_r\cdot a=p_ja.$$ The~order
idempotents $p_j$ are non-zero. Indeed, if $p_j=0$ for a~index $j=1,\ldots,m_{s_0}$ then
$ap_{j+1}=0$ and, hence, $p_{j+1}^{\rm d}ap_{j+1}=0$. Thus, $p_{j+1}=0$. Using the~elementary
induction, we conclude $p_j=0$ for all $j$, a~contradiction (we have shown that if
$p_1,\ldots,p_{m'}$ satisfies the~conditions of part~{\bf (b)} then $p_j\neq0$ for all
$j=1,\ldots,m'$ and, hence, $p_j\neq{\bf e}$ for all~$j$). Now we assume that there exists no such
index $s_0$. The~latter implies the~equality $\prod\limits_{k\in J_j}q_k=0$ for all $J\in{\cal
I}_s$ and all $s$, i.e., the~validity of the~relation~(\ref{12}) for all $J\in{\cal P}_{r-1}$.

Now, using the~identity~(\ref{10}), i.e., (\ref{12})~for $r=m$, and the~construction above and
taking the~finite number of steps, either we will construct the~required collection of order
idempotents $p_1,\ldots,p_{m'}$ or we will reduce the~problem to the~case when
the~identity~(\ref{10}) holds with $r=1$, i.e., to the~case of $q_1,\ldots,q_m=0$. The~latter is
impossible as $q_1\neq0$.

{\bf (b)}~$\Longrightarrow$~{\bf (a)} As was shown above, if $p_1,\ldots,p_{m'}$ satisfy
the~condition of part~{\bf (b)} then $p_j$ is not-trivial for all $j=1,\ldots,m'$. Now, using
the~elementary induction, we have
$$p_ja^{m'}=ap_{j+1}a^{m'-1}=a^2p_{j+2}a^{m'-2}=\ldots=a^{m'}p_{j+m'}=a^{m'}p_j$$
for an~arbitrary index~$j$, whence
$$p_ja^m=p_ja^{m'\frac{m}{m'}}=a^{m'}p_ja^{(m'-1)\frac{m}{m'}}=\ldots=a^mp_j.$$
Finally, $p_ja^mp_j^{\rm d}=0$ and, in~particular, $a^m$ is reducible.

{\bf (b)}~$\Longrightarrow$~{\bf (c)} We define the~element $d$ via the~formula
$d=\sum\limits_{j=1}^{m'}e^{i\frac{2\pi}{m'}j}p_j$, where the~order idempotents $p_1,\ldots,p_{m'}$
satisfy the~condition of part~{\bf (b)}. As is easy to see, the~element~$d$ is invertible and
$d^{-1}=\sum\limits_{r=1}^{m'}e^{-i\frac{2\pi}{m'}j}p_r$. We have the~equalities
$$e^{i\frac{2\pi}{m'}}dad^{-1}=
e^{i\frac{2\pi}{m'}}\Big(\sum\limits_{j=1}^{m'}e^{i\frac{2\pi}{m'}j}p_j\Big)a
\Big(\sum\limits_{r=1}^{m'}e^{-i\frac{2\pi}{m'}j}p_r\Big)=$$
$$=e^{i\frac{2\pi}{m'}}\sum\limits_{j,r=1}^{m'}e^{i\frac{2\pi}{m'}(j-r)}p_jap_r=
e^{i\frac{2\pi}{m'}}\sum\limits_{j=1}^{m'}e^{i\frac{2\pi}{m'}(j-(j+1))}ap_{j+1}=a.$$ Therefore,
$\sigma(a)=e^{i\frac{2\pi}{m'}}\sigma(dad^{-1})=e^{i\frac{2\pi}{m'}}\sigma(a)$. Consequently,
part~{\bf (b)} implies the~inva\-riance of the~spectrum $\sigma(a)$ under the~rotation on angle
$\frac{2\pi}{m'}$. In~particular, the~points $r(a)e^{i\frac{2\pi}{m'}j}$ belong to $\sigma_{\rm
per}(a)$ for $j=0,\ldots,m'-1$ as $r(a)\in\sigma(a)$. Let us show that these points are simple
poles of $R(\cdot, a)$. Using the~identity $\lambda-a=e^{i\frac{2\pi}{m'}}d(\lambda
e^{-i\frac{2\pi}{m'}}-a)d^{-1}$ for all $\lambda\in{\Bbb C}$, we obtain
$R(\lambda,a)=e^{-i\frac{2\pi}{m'}}dR(\lambda e^{-i\frac{2\pi}{m'}},a)d^{-1}$ for all
$\lambda\in\rho(a)$. Let $\lambda_0$ be a~pole of $R(\cdot,a)$ of order $k$. Then
$\lambda_0e^{i\frac{2\pi}{m'}}$ is an~isolated point of $\sigma(a)$ and we have the~equalities
$$\lim\limits_{\lambda\to\lambda_0e^{i\frac{2\pi}{m'}}}
(\lambda-\lambda_0e^{i\frac{2\pi}{m'}})^kR(\lambda,a)=
\lim\limits_{\lambda\to\lambda_0e^{i\frac{2\pi}{m'}}} (\lambda-\lambda_0e^{i\frac{2\pi}{m'}})^k
e^{-i\frac{2\pi}{m'}}dR(\lambda e^{-i\frac{2\pi}{m'}},a)d^{-1}=$$
$$=\lim\limits_{\lambda\to\lambda_0e^{i\frac{2\pi}{m'}}}
(\lambda e^{-i\frac{2\pi}{m'}}-\lambda_0)^ke^{i\frac{2\pi}{m'}k} e^{-i\frac{2\pi}{m'}}dR(\lambda
e^{-i\frac{2\pi}{m'}},a)d^{-1}=$$
$$=\lim\limits_{\mu\to\lambda_0}(\mu-\lambda_0)^ke^{i\frac{2\pi}{m'}(k-1)}dR(\mu,a)d^{-1}=0.$$
Thus, the point $\lambda_0e^{i\frac{2\pi}{m'}}$ is a~pole of $R(\cdot,a)$ of order which is not
greater than $k$. On the~other hand, $r(a)$ is a~simple pole of $R(\cdot,a)$ and so the~points
$r(a)e^{i\frac{2\pi}{m'}j}$, $j=0,\ldots,m'-1$, are also simple poles of $R(\cdot,a)$.

{\bf (c)}~$\Longrightarrow$~{\bf (a)} For $j=0,1,\ldots,m''-1$, we put
$\omega_j=e^{i\frac{2\pi}{m''}j}$. Using the~identity \mbox{\cite[p. 22]{BD}}
$\frac1{1-z^{m''}}=\frac1{m''}\sum\limits_{j=0}^{m''-1}\frac1{1-\omega_j^{-1}z}$ which is valid for
all complex numbers $z\notin\{\omega_0,\ldots,\omega_{m''-1}\}$, we have
$\frac1{\lambda^{m''}-z^{m''}}=
\frac1{m''\lambda^{m''-1}}\sum\limits_{j=0}^{m''-1}\frac{\omega_j}{\lambda\omega_j-z}$ for all
$\lambda,z\in{\Bbb C}$ satisfying $\frac{z}\lambda\notin\{\omega_0,\ldots,\omega_{m''-1}\}$.
Therefore, for $\lambda$ from a~sufficiently small punctured neighbourhood of $r(a)$, we obtain
\begin{equation}\label{15}
m''\lambda^{m''-1}R(\lambda^{m''},a^{m''})= \sum\limits_{j=0}^{m''-1}\omega_jR(\lambda\omega_j,a).
\end{equation}
In view of our condition and Lemma~\ref{lem10}, $r(a)^{m''}$ is a pole of $R(\cdot,a^{m''})$.
Assume that the~element $a^m$ is irreducible. Then the~element $a^{m''}$ is also irreducible.
Consequently, in view of Axiom~{\bf (A$_1$)}, $r(a)>0$, the~point $r(a)^m$ is a~simple pole of
the~resolvent $R(\cdot,a^{m''})$, and the~residue $(a^{m''})_{-1}\gg0$. The~point $r(a)$ is also
a~simple pole of the~function $\lambda\to R(\lambda^{m''},a^{m''})$. Indeed, $r(a)$ is an~isolated
singular point of this function and for every natural $k\in{\Bbb N}$ the~relation
\begin{equation}\label{16}
(\lambda-r(a))^kR(\lambda^{m''},a^{m''})=
\frac{(\lambda^{m''}-r(a)^{m''})^kR(\lambda^{m''},a^{m''})}
{(\lambda^{m''-1}+\lambda^{m''-2}r(a)+\ldots+r(a)^{m''-1})^k}
\end{equation}
holds. Therefore, for every $k>1$, we have
\begin{equation}\label{17}
\lim\limits_{\lambda\to r(a)}(\lambda-r(a))^kR(\lambda^{m''},a^{m''})=0.
\end{equation}
Let $k_j$ be the~order of the~pole of $R(\cdot,a)$ at the~point $r(a)\omega_j$, where
$j=0,1,\ldots,m''-1$. We claim equality $k_j=1$ (this fact is not assumed in part~{\bf (c)};
moreover, we mention that the~cone $A^+$ is not assumed to be normal). Put $l=\max\limits_{0\le
j\le m''-1}k_j$. Proceeding by contradiction, we suppose $l>1$. Therefore, using the~identities
(\ref{15}) and~(\ref{17}), we obtain
$$0=\sum\limits_{j=0}^{m''-1}\omega_j\lim\limits_{\lambda\to r(a)}
(\lambda-r(a))^lR(\lambda\omega_j,a)=$$
$$=\sum\limits_{j=0}^{m''-1}\lim\limits_{\lambda\to r(a)}
(\lambda\omega_j-r(a)\omega_j)^l\frac1{\omega_j^{l-1}}R(\lambda\omega_j,a)=
\sum\limits_{j=0}^{m''-1}\frac1{\omega_j^{l-1}}a_{-l,r(a)\omega_j}$$ and so
\begin{equation}\label{18}
\sum\limits_{j=0}^{m''-1}\frac1{\omega_j^{l-1}}a_{-l,r(a)\omega_j}=0.
\end{equation}
Taking into account the~integral representation of the~coefficients of the~Laurent series expansion
of the~resolvent $R(\cdot,a)$ around the~point $\omega_j$ and using the~functional calculus, for
arbitrary~$j_0$ satisfying $k_{j_0}=l$, we get
$a_{-1,r(a)\omega_{j_0}}a_{-l,r(a)\omega_{j_0}}=a_{-l,r(a)\omega_{j_0}}$ and
$a_{-1,r(a)\omega_{j_0}}a_{-l,r(a)\omega_j}=0$ for $j\neq j_0$. Now, from the~(\ref{18}), it
follows that
$$0=a_{-1,r(a)\omega_{j_0}}\sum\limits_{j=0}^{m''-1}\frac1{\omega_j^{l-1}}a_{-l,r(a)\omega_j}=
\frac1{\omega_{j_0}^{l-1}}a_{-l,r(a)\omega_{j_0}}.$$ Therefore, $a_{-l,r(a)\omega_{j_0}}=0$,
a~contradiction. Thus, $l=1$. Next, according to~(\ref{16}),
$$\lim\limits_{\lambda\to r(a)}(\lambda-r(a))R(\lambda^{m''},a^{m''})=
\frac{(a^{m''})_{-1}}{m''r(a)^{m''-1}}.$$ Using the~last equality and the~identity~(\ref{16}) once
more, we obtain
$$(a^{m''})_{-1}=\lim\limits_{\lambda\to r(a)}(\lambda-r(a))
\sum\limits_{j=0}^{m''-1}\omega_jR(\lambda\omega_j,a)=
\sum\limits_{j=0}^{m''-1}a_{-1,r(a)\omega_j}$$ and, hence,
$$(a^{m''})_{-1}a_{-1,r(a)\omega_j}=a_{-1,r(a)\omega_j}(a^{m''})_{-1}=a_{-1,r(a)\omega_j}.$$
In other  words, for all $j=0,1,\ldots,m''-1$ the idempotents $a_{-1,r(a)\omega_j}$ belong to
the~space $N(I-R_{(a^{m''})_{-1}})\cap N(I-L_{(a^{m''})_{-1}})$. In view of Axiom~{\bf (A$_2$)}
(see Proposition~\ref{prop9} and, in~particular, the~identity~(\ref{4})), this space is
one-dimensional. Therefore, taking into account the~condition $m''>1$, we get
$a_{-1,r(a)\omega_j}=0$ for all $j$. This contradiction establishes that $a^m$ is reducible.

Now we suppose that Axioms {\bf (A$_3$)} and~{\bf (A$_4$)} hold and show that for some $m_0\in{\Bbb
N}$, we have the~equality $\sigma_{\rm per}(a)\cap\sigma_{\rm j}(a;a_1)=r(a)H_{m_0}$. As is easy to
see, in~view of Axiom~{\bf (A$_4$)}, the~complexification $(A_{a_{-1}})_{\Bbb C}$ of the~order
ideal $A_{a_{-1}}$ generated by $a_{-1}$ is isomorphic onto the~complex linear subspace
$A_0=\{b+ic:b,c\in A_{a_{-1}}\}$ of $A$. On~the~other hand, in view of Axiom~{\bf (A$_3$)},
the~order ideal $A_{a_{-1}}$ is lattice isomorphic onto a~space $C(K)$ and, under this isomorphism,
$a_{-1}$~is mapped onto~${\rm 1}\!\!{\rm 1}_K$. In~particular, $A_{a_{-1}}$ is a~Riesz space and so
the~space $(C(K))_{\Bbb C}$ can be identified with~$A_0$, i.e., with the~complexification
$(A_{a_1})_{\Bbb C}$ of $A_{a_{-1}}$. If for a~number $\lambda\ge0$ and $b\in A_{a_{-1}}$,
the~inequalities $-\lambda a_{-1}\le b\le\lambda a_{-1}$ hold then, taking into account
the~identity $aa_{-1}=r(a)a_{-1}$, we have $-\lambda r(a)a_{-1}\le L_ab\le\lambda r(a)a_{-1}$.
The~latter implies the~\mbox{$L_a$-invariance} of $A_{a_{-1}}$ and the~relation
$\|L_ab\|_{a_{-1}}\le r(a)\|b\|_{a_{-1}}$. Whence $\|\widehat{L_a}\|_{B(A_{a_{-1}})}=r(a)$, where
the~positive operator $\widehat{L_a}$ is a~restriction of $L_a$ to $A_{a_{-1}}$. Thus,
$r(\widehat{L_a})=r(a)$. Since $a_{-1}$ is a~minimal idempotent (Axiom~{\bf (A$_2$)}) for arbitrary
$b\in A$, we find a~scalar $f(b)$ satisfying $a_{-1}ba_{-1}=f(b)a$. The~function $f$ on
$A_{a_{-1}}$ defined above is linear and positive and satisfies the~inequality
$|f(b)|\le\|b\|_{a_{-1}}$ for all $b\in A_{a_{-1}}$. In~particular, the~functional~$f$ is bounded.
Since $a_{-1}\gg0$, $f$ is strictly positive. For arbitrary $a,b\in A_{a_{-1}}$, we have
the~equalities
$$f(ab)a_{-1}=a_{-1}aba_{-1}=r(a)a_{-1}ba_{-1}=r(a)f(b)a_{-1}.$$
Whence for the~adjoint operator $\widehat{L_a}^*$, we get $(\widehat{L_a}^*f)(b)=f(ab)=r(a)f(b)$ or
$\widehat{L_a}^*f=r(a)f$.

Fix $\lambda_0\in\sigma_{\rm per}(a)\cap\sigma_{\rm j}(a;a_1)$. There exists an~element $x_0$ of
$A$ which is not nilpotent, satisfies the~relations $L_ax_0=R_ax_0=\lambda_0x_0$, and has
the~representation in the~form $x_0=x_0'+ix_0''$, where $x_0',x_0''\in A_{a_{-1}}$. In~particular,
$x_0\in(A_{a_{-1}})_{\Bbb C}$ and $\widehat{L_a}x_0=\lambda_0x_0$. For the~element $x_0$
the~modulus $|x_0|$ exists, belongs to $A_{a_{-1}}$, and is given by
$|x_0|=\sup\limits_{\varphi\in[0,2\pi]}{|(\cos{\varphi})x_0'+(\sin{\varphi})x_0'|}$. We can assume
$\||x_0|\|_{a_{-1}}=1$. For the~element $ax_0=ax_0'+iax_0''\in(A_{a_{-1}})_{\Bbb C}$ and for
the~modulus of this element in $(A_{a_{-1}})_{\Bbb C}$, we have
$$|ax_0|=\sup\limits_{\varphi\in[0,2\pi]}{|(\cos{\varphi})ax_0'+(\sin{\varphi})ax_0'|}\le
a\sup\limits_{\varphi\in[0,2\pi]}{|(\cos{\varphi})x_0'+(\sin{\varphi})x_0'|}=a|x_0|.$$ Whence
$r(a)|x_0|=|\lambda_0x_0|=|ax_0|\le a|x_0|$ or $0\le(a-r(a))|x_0|$. On the~other hand,
$0=a_{-1}(a-r(a))|x_0|$ and $a_{-1}\gg0$. Thus, $r(a)|x_0|=a|x_0|$; analogously,
$r(a)|x_0|=|x_0|a$. Therefore, for $\lambda>r(a)$, we have
$$R(\lambda,a)|x_0|=\frac1\lambda|x_0|+\frac1{\lambda^2}r(a)|x_0|+\ldots=
\frac{|x_0|}\lambda\Big(1+\frac{r(a)}\lambda+\Big(\frac{r(a)}\lambda\Big)^2+\ldots\Big)=
\frac{|x_0|}{\lambda_0-r(a)}$$ and, hence, $a_{-1}|x_0|=|x_0|$; analogously, $|x_0|a_{-1}=|x_0|$.
Consequently, we have the~equality $|x_0|=a_{-1}|x_0|=a_{-1}|x_0|a_{-1}=\mu a_{-1}$ for some
$\mu\in{\Bbb R}$. Since $\||x_0|\|_{a_{-1}}=1$, we infer $\mu=1$ and so $|x_0|=a_{-1}$. There
exists (see \mbox{\cite[Lemma 5.1(I)]{NS}}) an operator~$S$ on~$(A_{a_{-1}})_{\Bbb C}$ depending
upon~$x_0$ and satisfying the~identity $\widehat{L_a}=r(a)\lambda_0^{-1}S^{-1}\widehat{L_a}S$;
analogously, $\widehat{R_a}=r(a)\lambda_0^{-1}S^{-1}\widehat{R_a}S$  for the~same operator $S$,
where $\widehat{R_a}$ is a~restriction of~$R_a$ to $A_{a_{-1}}$. Now, if $\lambda'\in\sigma_{\rm
per}(a)\cap\sigma_{\rm j}(a;a_1)$ and $L_ax=R_ax=\lambda'x$ for $x\in(A_{a_{-1}})_{\Bbb C}$ then
$\widehat{L_a}Sx=\frac{\lambda_0\lambda'}{r(a)}Sx$ and
$\widehat{R_a}Sx=\frac{\lambda_0\lambda'}{r(a)}Sx$. Whence $\frac{\lambda_0\lambda'}{r(a)}$ belongs
to the~set $\sigma_{\rm per}(a)\cap\sigma_{\rm j}(a;a_1)$. Finally, since $\lambda_0$ and
$\lambda'$ are arbitrary, the set $\frac1{r(a)}\sigma_{\rm per}(a)\cap\sigma_{\rm j}(a;a_1)$ is
the~group of all $m_0^{\rm th}$~roots of unity for some $m_0\in{\Bbb N}$. We used Axiom~{\bf
(A$_4$)} where $b$ is an~algebraically strictly positive idempotent only. We didn't use
the~assumption that $x_0$ is not nilpotent (this assumption from the~definition of $\sigma_{\rm
j}(a;b)$ will be needed for the~validity of Lemma~\ref{lem15} and Corollary~\ref{cor16}).

The~proof of the~theorem is now complete.\hfill$\Box$
\smallskip

In the~proof of the~preceding theorem, using a~number $m$ of part~{\bf (a)}, a~number~$m'$
satisfying the~conditions of part~{\bf (b)} was found. A~number $m'$ is not uniquely
de\-ter\-mined. Indeed, for the~cyclic matrix $\left(\begin{smallmatrix} 0 & 1 & 0 & 0 \\ 0 & 0 & 1
& 0 \\ 0 & 0 & 0 & 1 \\ 1 & 0 & 0 & 0 \end{smallmatrix}\right)$ two situations $m'=2$ and $m'=4$
are possible (then the~required order idempotents of part~{\bf (b)} correspond, in~the~former, to
the~coordinates $\{1,3\}$ and $\{2,4\}$ and, in the~second case, to $\{1\},\ldots,\{4\}$). For
the~case of the~analogous cyclic $6\times6$ matrix three situations $m'=2$ ($\{1,3,5\}$ and
$\{2,4,6\}$), $m'=3$ ($\{1,4\}$, $\{2,5\}$, and $\{3,6\}$), and $m'=6$ ($\{1\}$,\ldots,$\{6\}$) are
possible.

As is easy to see from the~proof of Theorem~\ref{thm12}, if $p_1,\ldots,p_{m'}$ satisfying
the~conditions of part~{\bf (b)} exist then part~{\bf (c)} holds with $m''=m'$. In is not known if
the~converse is valid. I.e., does part~{\bf (b)} hold with $m''=m'$ if part~{\bf (c)} holds? In
other words, the~author does not know the~direct proof of the~implication {\bf
(c)}~$\Longrightarrow$~{\bf (b)} of the~preceding theorem. Nevertheless, in the~following
particular case such a~proof is possible. Let the~condition of part~{\bf (c)} of
Theorem~\ref{thm12} hold and let the~number $m''$ have the~representation in the~form of
the~product $m''=m_1\ldots m_n$, where~$n\in{\Bbb N}$, $m_k$ are prime numbers for all
$k=1,\ldots,n$, and $m_{k'}\neq m_{k''}$ for all $k'\neq k''$. Obviously, for arbitrary
$k=1,\ldots,n$ and all $j=0,1,\ldots,m_k-1$ the~points $r(a)e^{i\frac{2\pi}{m_k}j}$ are poles of
$R(\cdot,a)$ and so (see the~proof of the~implication {\bf (c)}~$\Longrightarrow$~{\bf (a)})
$a^{m_k}$ is reducible. In view of the~validity of the~implication {\bf (a)}~$\Longrightarrow$~{\bf
(a)} and since $m_k$ is a~prime number, there exists non-zero order idempotents
$p_{m_k,1},\ldots,p_{m_k,m_k}$ determining the~cyclic form of the~element $a$. For an~arbitrary
$n$-tuple $(j_1,\ldots,j_n)$, where $j_k=1,\ldots,m_k$ for all $k=1,\ldots,n$, we put
$p_{j_1,\ldots,j_n}=p_{m_1,j_1}\ldots p_{m_n,j_n}$. If $(j_1',\ldots,j_n')\neq(j_1'',\ldots,j_n'')$
then $p_{j_1',\ldots,j_n'}\perp p_{j_1'',\ldots,j_n''}$ and, hence, the~number of order idempotents
$p_{j_1,\ldots,j_n}$ is equal to~$m''$. Next,
$\sum\limits_{j_1=1}^{m_1}\ldots\sum\limits_{j_n=1}^{m_n}p_{j_1,\ldots,j_n}={\bf e}$,
$p_{j_1,\ldots,j_n}a=ap_{j_1-1,\ldots,j_n-1}$, and $p_{j_1-m'',\ldots,j_n-m''}=p_{j_1,\ldots,j_n}$.
If for some $r\in{\Bbb N}$ the~equality $p_{j_1-r,\ldots,j_n-r}=p_{j_1,\ldots,j_n}$ is valid then
for all $k=1,\ldots,n$ the~number $m_k$ is a~divisor of $r$. Therefore, $r=m_1\ldots m_kl=m''l$
with~$l\in{\Bbb N}$. Finally, the~collection $\{p_{j_1,\ldots,j_n}\}$ determines the~cyclic form of
$a$.
\smallskip

The~next result makes more precisely correlations between the~values $m$, $m'$, and~$m''$.

\begin{cor}\
Under the~assumptions of Theorem~{\rm \ref{thm12}}, the~following equalities hold
$$\min{\{m\in{\Bbb N}:a^m \ \text{is reducible}\}}=$$
$$=\min{\{m\in{\Bbb N}: \text{there exist} \ p_1,\ldots, p_m \ \text{in} \ {\bf OI}(A) \
\text{determining the~cyclic form of} \ a\}}$$
$$=\min{\{m\in{\Bbb N}:
r(a)e^{i\frac{2\pi}mj} \ \text{is a~pole of} \ R(\cdot,a) \ \text{for} \ j=0,1,\ldots,m-1\}};$$ if
$a^m$ is irreducible for all $m$ or such $p_1,\ldots, p_m$ does not exist then we suppose that
the~respective minimum is equal to one.
\end{cor}

We continue with the~following two auxiliary results.

\begin{lem}\label{lem14}\
Let $c\in A^+$, let $p_0,p_1,\ldots,p_n$, where $n\in{\Bbb N}$, be a~collection of
\mbox{$c$-invariant} order idempotents satisfying the~relations ${\bf e}=p_n\ge\ldots\ge p_1\ge
p_0$, and let $q_j=p_jp_{j-1}^{\rm d}$ for all $j=1,\ldots,n$. If for every $j$ the~elements
$c_{q_j}$ are nilpotent then the~element $c$ is also nilpotent.
\end{lem}

{\bf Proof.} For $n=1$ the~assertion is obvious. Let us assume $n>1$. We consider first the~case of
$c_{q_j}=0$ for all $j$ and will show the~identity $c^n=0$. To this end, using the~induction on
$k=1,\ldots,n$, we shall prove the~validity of the~equality $q_sc^kq_j=0$ for all $j=1,\ldots,n$
and $s=\max{\{1,j-k+1\}},\ldots,n$ which, for $k=n$, implies $q_sc^nq_j=0$ for all $s$ and $j$. In
view of the~relation $\sum\limits_{j=1}^nq_j={\bf e}$, we obtain $c^n=0$. If $j<s$ then $0\le
q_scq_j\le p_{s-1}^{\rm d}cp_j\le p_{j-1}^{\rm d}cp_j=0$. Taking into account our condition, we
have the~equality $q_scq_j=0$ for $j\le s$. Therefore, for $k=1$ our induction hypothesis is true
and the~identity $c=\sum\limits_{s=1}^{n-1}\sum\limits_{j=s+1}^nq_scq_j$ holds. Assume the~validity
of our assertion for some $k<m$. Since $n-k+1\ge\max{\{1,j-k+1\}}$ for all $j$, we obtain
$q_sc^kq_j=0$ for $s\ge n-k+1$ and this equality also holds for $j<s+k$. Therefore,
$c^k=\sum\limits_{s=1}^{n-k}\sum\limits_{j=s+k}^nq_sc^kq_j$. Evidently, $q_mc^{k+1}q_j=0$ for all
$j$ and $q_sq_j=0$ for $s\neq j$. Now, for arbitrary indexes $l$ and $t$ with $l\le m-1$, we have
$$q_lc^{k+1}q_t=q_lcc^kq_t=
\sum\limits_{j=l+1}^nq_lcq_j\cdot\sum\limits_{s=1}^{n-k}\sum\limits_{j=s+k}^nq_sc^kq_jq_t=
\sum\limits_{j=l+1}^nq_lcq_j\cdot\sum\limits_{s=1}^{t-k}q_sc^kq_t.$$ Consequently, if $t-k\le l$
then $q_lc^{k+1}q_t=0$, as desired.

In a~general case, we mention first the~validity of the~equality $(c^k)_{q_j}=(c_{q_j})^k$ for all
$k\in{\Bbb N}$ and $j=1,\ldots,n$. Indeed, the~case of $k=1$ is obvious. If the~required equality
is true for some~$k$ then
$$(c_q)^{k+1}=(c_q)^kc_q=(c^k)_qc_q=qc^kqcq=qc^{k+1}q-qc^kq^{\rm d}cq=(c^{k+1})_q.$$
Now we choose $k_j\in{\Bbb N}$ satisfying $c_{q_j}^k=0$ for all $j=1,\ldots,n$. Putting
$k=\max\limits_{1\le j\le n}{k_j}$, we have $(c^k)_{q_j}=(c_{q_j})^k=0$. As was showed above,
$c^{kn}=0$.\hfill$\Box$
\smallskip

For a~non-zero element $q\in{\bf OI}(A)$, we put $A(q)=\{a\in A:qaq=a\}$. Under the~linear
operations, the~multiplication, the~norm, and the~order induced by $A$, the~linear space $A(q)$ is
an~ordered Banach algebra with unit $q$ being a~closed subalgebra of $A$. We mention at once
the~following properties of $A(q)$~\cite{Al2}: {\bf (a)} For the~order interval $[0,q]$ in $A$
the~identity ${\bf OI}(A(q))={\bf OI}(A)\cap[0,q]$ is valid; {\bf (b)} An~element $z\in A(q)$ is
irreducible with respect $q$ if and only if $z$ is irreducible in $A(q)$; {\bf (c)} The~identity
$(A(q))_{\rm n}=A_{\rm n}\cap A(q)$ holds; {\bf (d)} If $A$ has a~disjunctive product then $A(q)$
also has a~disjunctive product. From this remarks, it follows that if $A$ satisfies Axiom~{\bf
(A$_j$)} for some $j=1,3,4$ then $A(q)$ also satisfies this axiom.

\begin{lem}\label{lem15}\
Let $A$ satisfy Axiom~{\bf (A$_1$)}. Let an~element $a\in A_{\rm n}$ such that $r(a)>0$ is a~pole
of $R(\cdot,a)$ and the~residue $a_{-1}$ possesses a~modulus $|a_{-1}|$. Let the~collection of
order idempotents $p_0,p_1,\ldots,p_n$ determine the~Frobenius normal form of the~element~$a$. Then
the~inclusions
\begin{eqnarray}\label{19}
& \sigma_{\rm per}(a;A)\cap\sigma_{\rm j}(a;|a_{-1}|;A)\subseteq &\nonumber\\
& \subseteq \bigcup\{\sigma_{\rm per}(a_{q_j};A(q_j))\cap\sigma_{\rm
j}(a_{q_j};(a_{q_j})_{-1};A(q_j)): r(a_{q_j})=r(a)\}\subseteq\sigma_{\rm per}(a;A) &
\end{eqnarray}
hold with $q_j=p_jp_{j-1}^{\rm d}$ for $j=1,\ldots,n$.
\end{lem}

{\bf Proof.} We denote the~middle part of~(\ref{19}) by ${\cal S}$. By Corollary~\ref{cor5}, to
check the~left inclusion in~(\ref{19}), it suffices to establish the~relation
$$\bigcup\{\sigma_{\rm per}(a_{q_j};A)\cap\sigma_{\rm j}(a;|a_{-1}|;A): r(a_{q_j})=r(a)\}\subseteq
{\cal S}.$$ To this end, let $\lambda\in\sigma_{\rm j}(a;|a_{-1}|;A)$ and let $|\lambda|=r(a)$.
Then there exists an~element $c\in A$ which is not nilpotent, satisfies $ac=ca=\lambda c$, and has
a~representation in the~form $c=c'+ic''$ with $c',c''\in A_{|a_{-1}|}$. For arbitrary
$j=0,1,\ldots,n$, the~order idempotent~$p_j$ is \mbox{$|a_{-1}|$-invariant} and, hence, is
\mbox{$c$-invariant}. Therefore, for $j=1,\ldots,n$, we have
$$q_j^{\rm d}cq_j=(p_j^{\rm d}+p_{j-1})cp_jp_{j-1}^{\rm d}=p_{j-1}cq_j;$$
analogously, $q_jaq_j^{\rm d}=q_jcp_j^{\rm d}$. Thus, $q_jaq_j^{\rm d}cq_j=0$. Consequently,
$$\lambda c_{q_j}=q_jacq_j=q_jaq_jcq_j+q_jaq_j^{\rm d}cq_j=a_{q_j}c_{q_j};$$
analogously, $\lambda c_{q_j}=c_{q_j}a_{q_j}$. In view of Lemma~\ref{lem14}, there exists an~index
$j_0$ such that $c_{q_{j_0}}$~is not nilpotent. In~particular, $c_{q_{j_0}}\neq0$ and, hence,
$\lambda\in\sigma(a_{q_{j_0}})$. On the~other hand, Lemma~\ref{lem3} yields $r(a_{q_{j_0}})\le
r(a)$. Therefore, $r(a_{q_{j_0}})=r(a)$, $a_{q_{j_0}}$ is irreducible with respect~$q_{j_0}$, and
$\lambda\in\sigma_{\rm per}(a_{q_{j_0}};A(q_{j_0}))$. Then \cite{Al2} $r(a)$ is a~simple pole of
$R(\cdot,a_{q_{j_0}})$ and for $\lambda$ from a~sufficient small punctured neighbourhood of $r(a)$
the~identity
$$R(\lambda,a_{q_{j_0}})=
\frac1\lambda(p_{j_0}^{\rm d}+p_{j_0-1})+q_{j_0}R(\lambda,a)q_{j_0}$$ holds. Thus,
$0\le(a_{q_{j_0}})_{-1}=(a_{-1})_{q_{j_0}}$. The~last equality implies
$(c')_{q_{j_0}},(c'')_{q_{j_0}}\in A_{(a_{q_{j_0}})_{-1}}$ and, hence, $\lambda\in\sigma_{\rm
j}(a_{q_{j_0}};(a_{q_{j_0}})_{-1};A(q_{j_0}))$.

To check the~right inclusion in~(\ref{19}), mention the~next inclusion $\sigma_{\rm
per}(a_q;A(q))\subseteq\sigma_{\rm per}(a_q;A)$ for arbitrary $q\in{\bf OI}(A)$. In fact, if
$(\lambda-a_q)b=b(\lambda-a_q)={\bf e}$ for an~element $b\in A$ then
$(\lambda-a_q)b_q=b_q(\lambda-a_q)=q$. Now it only remains to recall Lemma~\ref{lem3}.\hfill$\Box$
\smallskip

Theorem~\ref{thm12} and the~preceding lemma imply the~next consequence which charac\-terizes
the~peripheral spectrum of an~arbitrary positive element of $A$.

\begin{cor}\label{cor16}\
Let the~assumptions of Lemma~{\rm \ref{lem15}} satisfy, let $A$ satisfy Axioms {\bf (A$_3$)}
and~{\bf (A$_4$)}, and let the~algebra $A(q)$ satisfy Axiom~{\bf (A$_2$)} for every non-zero
$q\in{\bf OI}(A)$. Then the~inclusions
$$\sigma_{\rm per}(a;A)\cap\sigma_{\rm j}(a;|a_{-1}|;A)\subseteq
r(a)\bigcup\limits_{s=1}^nH_{m_s}\subseteq\sigma_{\rm per}(a;A)$$ hold with some
$m_1,\ldots,m_n\in{\Bbb N}$.
\end{cor}

We close this section with the~next assertion about elements having the~cyclic form.

\begin{prop}\
Let an~element $a\in A$ have the~cyclic form. Suppose that there exists at least one pole of
the~resolvent $R(\cdot,a)$ among points of $\sigma_{\rm per}(a)$. Then the~decomposition
$a^m=\sum\limits_{j=1}^mb_j$ holds, where $1<m\in{\Bbb N}$, $b_{j'}b_{j''}=0$ for $j'\neq j''$, and
$r(b_j)=r(a)^m$ for all $j=1,\ldots,m$. If, in~addition, $a\in A^+$ then there exists such
a~decomposition that $b_{j'}\wedge b_{j''}=0$ for $j'\neq j''$.
\end{prop}

{\bf Proof.} Using our condition, we find order idempotents $p_1,\ldots,p_m$ of $A$, where $m>1$,
satisfying $p_ja=ap_{j+1}$ for $j=1,\ldots,m$. Evidently, $p_ja^m=ap_{j+1}a^{m-1}=\ldots=a^mp_j$,
whence $p_ja^mp_j=a^mp_j=p_ja^m$ for all $j$ and $p_{j'}a^mp_{j''}=0$ for $j'\neq j''$. Therefore,
putting $b_j=p_ja^mp_j$, we obtain $a^m=\sum\limits_{j=1}^mb_j$. Fix an~index $j$. Let $\lambda_0$
be a~pole of~$R(\cdot,a)$ and let $|\lambda_0|=r(a)$. Since the~equality $p_jR(\lambda,a)=0$ is
impossible $\lambda$ for sufficiently close to $\lambda_0$, $\lambda\neq\lambda_0$, there exists
a~number $s\in{\Bbb N}$ satisfying $p_ja_{\lambda_0,s}\neq0$ and $p_ja_{\lambda_0,t}=0$ for $t<s$.
Keeping the~identity $aa_{\lambda_0,s}=r(a)a_{\lambda_0,s}+a_{\lambda_0,s-1}$ in mind and using
the~elementary induction, it is not difficult to check the~validity of the~relation
$a^na_{\lambda_0,s}=\sum\limits_{t=0}^n\lambda_0^tC_n^ta_{\lambda_0,s-(n-t)}$ for all $n\in{\Bbb
N}$, where $C_n^t$ are binomial coefficients. The~latter implies the~equality
$p_ja^ma_{\lambda_0,s}=\lambda_0^mp_ja_{\lambda_0,s}$ or
$b_ja_{\lambda_0,s}=\lambda_0^ma_{\lambda_0,s}$. Thus, $\lambda_0^m\in\sigma(b_j)$. Taking into
account the~inequality $r(b_j)\le r(a)^m$, we have $r(b_j)=r(a)^m$. In view of the~definition of
$b_j$, the~last assertion is clear.\hfill$\Box$

\section{The~Lotz-Schaefer axiom}

In the~spectral theory of positive operators the~next Lotz-Schaefer theorem is one of the~most
significant results (see \mbox{\cite[pp. 351-352]{Sch}}). In some situations, e.g., in~the~case of
irreducible operators (see~\cite{Al}), this theorem allows the~study of points of the~peripheral
spectrum of an~operator to reduce to the~case of poles only.

\begin{thm}\
Let $T$ be a~positive operator on a~Banach lattice $E$ and let $r(T)$ be a~finite-rank pole of
$R(\cdot,T)$. Then the~peripheral spectrum $\sigma_{\rm per}(T)$ of $T$ consists entirely of poles
of $R(\cdot,T)$.
\end{thm}

As was mentioned above (see Axiom~{\bf (A$_2$)} and the~remarks before this axiom), algebraically
strictly positive projections in an~ordered Banach algebra $A$ can be considered as
a~generalization of rank-one operators. Moreover, the~residue $a_{-1}$ of the~resolvent
$R(\cdot,a)$ of the~irreducible element $a$ at the~point $r(a)$ satisfies this condition (see
Theorem~\ref{thm6}{\bf (c)}). An~element $b\in A$ is said to be {\it relatively algebraically
strictly positive} whenever there exists a~non-zero order idempotent~$q$ of $A$ such that $b\in
A(q)$ and $q_1aq_2>0$ for all $0<q_1,q_2\in{\bf OI}(A)\cap[0,q]$; in this case, we write $b\gg_q0$.
Obviously, $b\gg0$ if and only if $b\gg_{\bf e}0$. As is easy to see, if $P$ is a~non-zero order
continuous projection on a~Dedekind complete Banach lattice $E$ being relatively algebraically
strictly positive element in the~ordered Banach algebra $B(E)$ then $\dim{R(P)}=1$. On~the~other
hand, the~collection $F(E)$ of finite-rank operators on $E$ is an~algebraic ideal of $B(E)$ and if
a~non-zero operator $T\in F(E)$ then the~algebraic ideal generated by~$T$ in $B(E)$ coincides
with~$F(E)$. Keeping these remarks in mind, we define the~set ${\cal F}(A)$ of {\it finite-rank
elements} of an~ordered Banach algebras $A$ as the~(two-sided) algebraic ideal generated by the~set
of all relatively algebraically strictly positive order continuous projectors of $A$; if such
projectors do not exist, we put ${\cal F}(A)=\emptyset$. If $a\in A(q)$ is a~relatively
algebraically strictly positive element of $A(q)$ then $a$ is such an~element of $A$, and, hence,
${\cal F}(A(q))\subseteq{\cal F}(A)$. Now, axiomatizing the~respective theorem, we can introduce
the~Lotz-Schaefer axiom in the~following manner:
\begin{description}
\item[{\bf (A$_{\rm LS}$)}] If $a$ is a~positive element of $A$, $r(a)$ is a~pole of $R(\cdot,a)$,
and the~residue $a_{-1}$ is a~finite-rank element then the~peripheral spectrum $\sigma_{\rm
per}(a)$ of $a$ consists entirely of poles of $R(\cdot,a)$.
\end{description}

If $E$ is a~Dedekind complete Banach lattice admitting a~weak order unit $x>0$ and a~strictly
positive order continuous functional $f$ then the~projection $\frac1{f(x)}f\otimes x\gg0$ in
the~ordered Banach algebra $B(E)$ and, hence, ${\cal F}(B(E))=F(E)$. Therefore, the~Lotz-Schaefer
theorem implies the~validity of Axiom~{\bf (A$_{\rm LS}$)} for a~wide class of ordered Banach
algebras of the~form $B(E)$. If we want the~validity of a~similar axiom for a~wider class of
algebras of the~form $B(E)$ then we must introduce the~next weaker axiom:
\begin{description}
\item[{\bf (A$_{\rm LS}'$)}] If $a$ is a~positive element of $A$, $r(a)$ is a~pole of $R(\cdot,a)$,
the~residue $a_{-1}$ is an~order continuous element, and $a_{-1}\gg0$ then the~peripheral spectrum
$\sigma_{\rm per}(a)$ of $a$ consists entirely of poles of $R(\cdot,a)$.
\end{description}

Now if $E$ is an~arbitrary Dedekind complete Banach lattice then the~ordered Banach algebra~$B(E)$
satisfies Axiom~{\bf (A$_{\rm LS}'$)}. Moreover, the~last two conditions of this axiom
automatically hold when the~element $a$ is irreducible (see Theorem~\ref{thm6}{\bf (c)}). However,
as the~results below show (see, e.g., Theorem~\ref{thm23}), in the~case of an~arbitrary positive
elements Axiom~{\bf (A$_{\rm LS}$)} is more useful than Axiom~{\bf (A$_{\rm LS}'$)}.

Nevertheless, there are other cases when an~algebra $A$ need not satisfy Axiom~{\bf (A$_{\rm
LS}$)}. If ${\bf OI}(A)=\{0,{\bf e}\}$ then ${\bf e}\gg0$ and, hence, ${\cal F}(A)=A$. Next, if
$A=L_\infty(\Omega,\mu)$, with $\mu$~a~\mbox{$\sigma$-finite} diffuse measure on
a~$\sigma$-algebra, then there exist no algebraically relatively strictly positive elements in $A$
and, hence, ${\cal F}(A)=\{0\}$.

If we will assume that the~residue~$a_{-1}$ is a~minimal idempotent of $A$ instead of
$a_{-1}\in{\cal F}(A)$ then Axiom~{\bf (A$_{\rm LS}$)} is valid for the~algebra $B(E)$, where $E$
is an~arbitrary Banach lattice. Unfortunately, the~residue $a_{-1}$ of the~resolvent $R(\cdot,a)$
of an~irreducible element $a$ which belongs to an~ordered Banach algebra $A$ satisfying Axiom~{\bf
(A$_1$)} need not be a~minimal idempotent (see~\cite{Al2}). Therefore, under this assumption,
Axiom~{\bf (A$_{\rm LS}$)} is employed only for a~narrow class of ordered Banach algebras.
Moreover, as the~next example shows, Axiom~{\bf (A$_{\rm LS}'$)} in its present form and Axiom~{\bf
(A$_{\rm LS}$)} with the~assumption about the~minimality of $a_{-1}$ do not hold in general for
an~arbitrary ordered Banach algebra $A$. This example also shows that {\it the~peripheral spectrum
$\sigma_{\rm per}(a)$ of an~irreducible element $a$ need not be cyclic} while $r(a)$ is a~pole of
$R(\cdot,a)$ and the~residue $a_{-1}\gg0$.

\begin{exm}\
{\rm Consider the~space $\ell_\infty$ of all bounded sequences $x=(x_1,x_2,\ldots)$. Under
the~natural algebraic operations, multiplication, and $\sup$-norm, this space is a~commutative
Banach algebra with unit ${\bf e}=(1,1,\ldots)$. Fix an~arbitrary number $\lambda_0\in{\Bbb C}$ and
a~sequence $\{z_n\}$ in ${\Bbb C}$ satisfying $|z_n|=1$, $z_n\neq\lambda_0$ for all $n$,
$\lambda_0\neq1$ , $z_1=1$, $z_n\neq z_m$ for all $n\neq m$, and $z_n\to\lambda_0$ as $n\to\infty$.
Define the~sequence $z\in\ell_\infty$ by $z=(z_1,z_2,\ldots)$ and consider the~algebraic wedge
$K_0$ generated by ${\bf e}$ and~$z$, i.e.,
$$K_0=\Big\{\sum\limits_{j=0}^n\alpha_jz^j:
n\in{\Bbb N} \ \text{and} \ \alpha_j\in{\Bbb R}^+ \ \text{for all} \ j=0,1,\ldots,n\Big\}.$$ We
claim that the~closure $\overline{K_0}$ of $K_0$ is a~normal cone. Indeed, for arbitrary
$y=\sum\limits_{j=0}^n\alpha_jz^j$ with $\alpha_j\in{\Bbb R}$ the~inequalities
$\Big|\sum\limits_{j=0}^n\alpha_j\Big|\le\|y\|_{\ell_\infty}\le\sum\limits_{j=0}^n|\alpha_j|$ holds
as $z_1=1$. Therefore, if $\alpha_j\ge0$ for $j=0,1,\ldots,n$ then
\begin{equation}\label{20}
\|y\|_{\ell_\infty}=\sum\limits_{j=0}^n\alpha_j.
\end{equation}
Let $x\in\ell_\infty$ and let $\pm x\in\overline{K_0}$. There exist two sequences $\{x_n\}$ and
$\{y_n\}$ in $K_0$ satisfying $x_n\to x$ and $y_n\to-x$ and, hence, $x_n+y_n\to0$ as $n\to\infty$.
Let $x_n=\sum\limits_{j=0}^{k_n}\alpha_{nj}z^j$ and $y_n=\sum\limits_{j=0}^{k_n}\beta_{nj}z^j$,
where $\alpha_{nj},\beta_{nj}\ge0$ for all $n$ and $j=1,\ldots,k_n$. Then
$\sum\limits_{j=0}^{k_n}(\alpha_{nj}+\beta_{nj})z^j\to0$ and so
$0\le\sum\limits_{j=0}^{k_n}\alpha_{nj}\le\sum\limits_{j=0}^{k_n}(\alpha_{nj}+\beta_{nj})\to0$ as
$n\to\infty$. Taking into account the~identity~(\ref{20}), we have
$\|x_n\|_{\ell_\infty}=\sum\limits_{j=0}^n\alpha_{nj}\to0$. Therefore, $x=0$. Thus,
$K=\overline{K_0}$ is a~cone and, under the~order induced by $K$, $\ell_\infty$ is an~ordered
Banach algebra. In view of the~relation $z_n\neq z_m$ with $n\neq m$, the~system $\{{\bf
e},z,z^2,\ldots\}$ is linearly independent and, in~particular, every element $w\in K_0$ has
a~unique representation in the~form $w=\sum\limits_{j=0}^n\omega_jz^j$ with $\omega_j\ge0$. Now if
$u=\sum\limits_{j=0}^m\mu_jz^j\in K_0$ and $0\le_{K_0}w\le_{K_0}u$ then
$\|w\|_{\ell_\infty}=\sum\limits_{j=0}^n\omega_j\le\sum\limits_{j=0}^n\mu_j=\|u\|_{\ell_\infty}$.
Thus, $K_0$ is a~normal cone and, hence \mbox{\cite[p. 81, Exercise 9]{AlT}}, $K$ is normal. We
mention at once that the~relation $\overline{K-K}\neq\ell_\infty$ holds as the~Banach space
$\ell_\infty$ is not separable; in~particular, $K$~is not generating. Obviously, $r(z)=1$ and
$R(\lambda,z)=(\frac1{\lambda-z_1},\frac1{\lambda-z_2},\ldots)$ for all
$\lambda\notin\sigma(z)=\overline{\{z_1,z_2,\ldots\}}$. In~particular, $\xi=1$ is an~isolated point
of $\sigma(z)$. For $\lambda$ close to this point and $\lambda\neq1$ the~inequality
$|\lambda-1|\le|\lambda-z_n|$ holds for all $n$. Whence
$\Big|\frac{(\lambda-1)^2}{\lambda-z_n}\Big|\le|\lambda-1|$ and so
$(\lambda-1)^2\|R(\lambda,z)\|_{\ell_\infty}\le|\lambda-1|\to0$ as~$\lambda\to1$. Thus, $\xi=1$ is
a~simple pole of $R(\cdot,z)$. Let $r=\inf\limits_{n>1}{|1-z_n|}>0$. For arbitrary
$\epsilon$~and~$\lambda$ satisfying $0<\epsilon<r$ and $|\lambda-1|<\epsilon$, we have
$|\lambda-z_n|\ge|1-z_n|-|\lambda-1|\ge r-\epsilon$ and, hence,
$\Big|\frac{\lambda-1}{\lambda-z_n}\Big|\le\frac{\epsilon}{r-\epsilon}\to0$ as $\epsilon\to0$.
Consequently, $\lim\limits_{\lambda\to1}\Big|\frac{\lambda-1}{\lambda-z_n}\Big|=0$ uniformly in
$n=2,3,\ldots$ and so the~residue $R(\cdot,z)$ at the~point $\xi=1$ is equal to ${\bf
e_1}=(1,0,0,\ldots)$. If we can verify the~identity ${\bf OI}(\ell_\infty)=\{0,{\bf e}\}$ (under
the~order induced by $K$) then this means the~validity of the~relation ${\bf e_1}\gg0$ while
the~point $\lambda_0$ is not an~isolated point of the~spectrum $\sigma(z)$. To this~end, let
the~sequence $p\in{\bf OI}(\ell_\infty)$. Then $p$ is a~characteristic function $\chi_A$ of
a~subset $A$ of ${\Bbb N}$ which has a~representation in the~form
$p=\sum\limits_{j=0}^n\theta_jz^j$ with $\theta_j\ge0$. If $1\notin A$ then $(\chi_A)_1=0$, whence
$\sum\limits_{j=0}^n\theta_j=0$ or $p=0$. If $1\in A$ then $1\notin{\Bbb N}\setminus A$, whence
$p={\bf e}$, as required. Obviously, the~residue ${\bf e_1}$ is a~minimal idempotent of the~algebra
$\ell_\infty$. Moreover, as is easy to see, we can choose the~sequence $\{z_n\}$ such that
the~peripheral spectrum $\sigma_{\rm per}(z)$ of the~irreducible element~$z$ is not
cyclic.}\hfill$\Box$
\end{exm}

Now, assuming the~Lotz-Schaefer axiom~{\bf (A$_{\rm LS}$)}, we continue the~study of
the~perip\-heral spectrum and, in~particular, will obtain some consequences of Theorem~\ref{thm12}.
\smallskip

Before, we discuss the~following property of the~algebra $B(E)$, where $E$ is a~complex Banach
lattice being the~complexification of the~real Banach lattice $E_{\Bbb R}$. Let
$\{S_n\}$~and~$\{T_n\}$ be two sequences in the~space $B_{\rm r}(E)=B_{\rm r}(E_{\Bbb R})$ of all
regular operators on $E$ such that the~sequence $\{S_n+iT_n\}$ converges in $B(E)$. Then
$\{S_n\}$~and~$\{T_n\}$ are also convergent. Indeed, fix a~number $\epsilon>0$ and find an~index
$k\in{\Bbb N}$ satisfying $\|S_n-S_m+i(T_n-T_m)\|_{B(E)}<\epsilon$ for all $n,m\ge k$. Using
the~condition $S_n\in B_{\rm r}(E)$ and the~inequality $\|y+iz\|_E\ge\max{\{\|y\|_E,\|z\|_E\}}$ for
all $y,z\in E_{\Bbb R}$, we have
$$\|(S_n-S_m)x\|_E\le\|(S_n-S_m)x+i(T_n-T_m)x\|_E<\epsilon$$
for an~arbitrary element $x\in E_{\Bbb R}$ with $\|x\|_E=\|x\|_{E_{\Bbb R}}=1$. Whence
$\|S_n-S_m\|_{B(E_{\Bbb R})}<\epsilon$ and so $\|S_n-S_m\|_{B(E)}<2\epsilon$. Thus, the~sequence
$\{S_n\}$ is convergent; the~case of $\{T_n\}$ is analogous. We axiomatize this property and make
the~next assumption:
\begin{description}
\item[{\bf (A$_5$)}] The~convergence of the~sequence $\{b_n+ic_n\}$, where $\{b_n\}$ and
$\{c_n\}$ are two sequences in~$A_{\rm r}$, implies the~convergence of $\{b_n\}$ and $\{c_n\}$.
\end{description}

Evidently, if an~ordered Banach algebra $A$ satisfies Axiom~{\bf (A$_5$)} then the~algeb\-ra~$A(q)$
also satisfies this axiom for every non-zero order idempotent $q$ of $A$. Next, as the~example of
the~ordered Banach algebra $C^1[a,b]$ of all complex functions $x$ represented in the~form
$x=x_1+ix_2$, where the~functions $x_1,x_2:[a,b]\to{\Bbb R}$ are continuously differentiable, under
the~natural algebraic operations, multiplication, order, and norm
$\|x\|_{C^1[a,b]}=\max\limits_{t\in[a,b]}|x(t)|+\max\limits_{t\in[a,b]}|{\dot x}(t)|$, shows,
Axiom~{\bf (A$_5$)} does not imply the~normality of a~cone~$A^+$.

\begin{thm}\label{thm20}\
Let an~ordered Banach algebra $A$ satisfy Axioms~{\bf (A$_1$)}-{\bf (A$_5$)} and~{\bf (A$_{\rm
LS}'$)}. Let an~element $a$ of $A$ be non-zero, order continuous, and irreducible. Let the~point
$r(a)$ be a~pole of $R(\cdot,a)$. Then the~peripheral spectrum $\sigma_{\rm per}(a)$ consists
entirely of poles of $R(\cdot,a)$ and has the~form $\sigma_{\rm per}(a)=r(a)H_m$ for some
$m\in{\Bbb N}$.
\end{thm}

{\bf Proof.} In view of Theorem~\ref{thm12}, it suffices to establish the~inclusion $\sigma_{\rm
per}(a)\subseteq\sigma_{\rm j}(a;a_{-1})$. To this end, let $\lambda_0\in\sigma_{\rm per}(a)$.
Since $A$ satisfies Axiom~{\bf (A$_{\rm LS}'$)}, $\lambda_0$ is a~pole of $R(\cdot,a)$ of order
$k$. The~relations $aa_{\lambda_0,-k}=a_{\lambda_0,-k}a=\lambda_0a$ hold. We claim that $k=1$ and
the~residue $a_{\lambda_0,-1}$ has the~representation in the~form $a_{\lambda_0,-1}=z'+iz''$, where
the~elements $z',z''\in A_{a_{-1}}$. The~last two assertions imply the~inclusion
$\lambda_0\in\sigma_{\rm j}(a;a_{-1})$. To check them, let
$\lambda_0=r(a)(\cos{\varphi}+i\sin{\varphi})$ with $\varphi\in[0,2\pi)$. We have the~equalities
$$a_{\lambda_0,-k}=\lim\limits_{\lambda\to\lambda_0}(\lambda-\lambda_0)^kR(\lambda,a)=
\lim\limits_{t\downarrow1}(t\lambda_0-\lambda_0)^kR(t\lambda_0,a)=
\lim\limits_{t\downarrow1}(t-1)^k\lambda_0^k\sum\limits_{j=0}^\infty\frac1{(t\lambda_0)^{j+1}}a^j
=$$
$$=\lim\limits_{t\downarrow1}(t-1)^k\sum\limits_{j=0}^\infty
\frac{\cos{((k-j+1)\varphi)}+i\sin{((k-j+1)\varphi)}}{t^{j+1}}a^j.$$ Fix an~arbitrary sequence
$\{t_n\}$ in ${\Bbb R}$ satisfying $t_n\downarrow1$ and put
$$b_n=(t_n-1)^k\sum\limits_{j=0}^\infty\frac{\cos{((k-j+1)\varphi)}}{t_n^{j+1}}a^j \ \ {\rm and} \ \
c_n=(t_n-1)^k\sum\limits_{j=0}^\infty\frac{\sin{((k-j+1)\varphi)}}{t_n^{j+1}}a^j.$$ Obviously,
$b_n+ic_n\to a_{\lambda_0,-k}$ as $n\to\infty$ and $\pm
b_n\le(t_n-1)^k\sum\limits_{j=0}^\infty\frac1{t_n^{j+1}}a^j\to a_{-k}$; analogously, for $\{c_n\}$.
Taking into account Axiom~{\bf (A$_5$)}, we conclude the~convergence of $\{b_n\}$ and $\{c_n\}$ and
the~validity of the~representation $a_{\lambda_0,-k}=x'+ix''$, where $-a_{-k}\le x',x''\le a_{-k}$.
If $k>1$ then $a_{-k}=0$ and, hence, $x',x''=0$ or $a_{\lambda_0,-k}=0$, a~contradiction. Thus,
$k=1$, as required.\hfill$\Box$

\begin{lem}\label{lem21}\
If $A$ satisfies Axiom~{\bf (A$_{\rm LS}$)} then $A(q)$ also satisfies this axiom for every
non-zero $q\in{\bf OI}(A)$.
\end{lem}

{\bf Proof.} Let $0\le a\in A(q)$, let $r(a)$ be a~pole of $R_q(\cdot,a)$, where $R_q(\cdot,a)$ is
a~resolvent of $a$ in~$A(q)$, and let the~residue $a_{-1}$ of $R_q(\cdot,a)$ at the~point $r(a)$
satisfy the~condition $a_{-1}\in{\cal F}(A(q))$. We can assume $r(a)>0$. The~inclusions~\cite{Al2}
\begin{equation}\label{21}
\rho_\infty(a;A)\subseteq\rho(a;A(q)) \ \ {\rm and} \ \
\rho(a;A(q))\setminus\{0\}\subseteq\rho(a;A)
\end{equation}
hold, where $\rho_\infty(a;A)$ is the~unbounded connected component of $\rho(a;A)$. Moreover, for
arbitrary $\lambda\in\rho(a;A(q))\setminus\{0\}$, we have the~identity~\cite{Al2}
\begin{equation}\label{22}
R(\lambda,a)=R_q(\lambda,a)+\frac1\lambda q^{\rm d}.
\end{equation}
Thus, $r(a)$ is a~pole of $R(\cdot,a)$ and the~residues of $R(\cdot,a)$ and $R_q(\cdot,a)$ at
$r(a)$ coincide. In view of the~inclusion ${\cal F}(A(q))\subseteq{\cal F}(A)$, $a_{-1}\in{\cal
F}(A)$. If $\lambda\in\sigma_{\rm per}(a;A(q))$ then $\lambda\in\sigma_{\rm per}(a;A)$ and, hence,
$\lambda$ is a~pole of $R(\cdot,a)$. Taking into account the~first inclusion of~(\ref{21}), we
infer that $\lambda$ is an~isolated point of $\sigma(a;A(q))$. Consequently, in view of~(\ref{22}),
$\lambda$ is a~pole of $R_q(\cdot,a)$.\hfill$\Box$
\smallskip

The~next result which follows easily from Corollary~\ref{cor5}, Theorem~\ref{thm20},
the~inclu\-sions~(\ref{21}), and the~preceding lemma, characterizes the~peripheral spectrum of
a~wide class of positive elements.

\begin{cor}\
Let $A$ satisfy Axioms {\bf (A$_1$)}, {\bf (A$_3$)}-{\bf (A$_5$)}, and~{\bf (A$_{\rm LS}$)}, let
the~algebra $A(q)$ satisfy Axiom~{\bf (A$_2$)} for every non-zero $q\in{\bf OI}(A)$, and let
an~element $a\in A_{\rm n}$ such that it has the~Frobenius normal form and $r(a)$ is a~pole of
$R(\cdot,a)$. Then the~identity $\sigma_{\rm per}(a)=r(a)\bigcup\limits_{s=1}^nH_{m_s}$ holds with
some $m_1,\ldots,m_n\in{\Bbb N}$.
\end{cor}

As the~next theorem shows, Axiom~{\bf (A$_{\rm LS}$)} can be employed to positive elements having
the~Frobenius normal form. Moreover, this result, Theorem~\ref{thm20}, and the~pre\-ceding
corollary illustrate the~importance of the~notion of finite-rank element as we define it above. It
also shows that our definition in the~abstract case is the~right one to use.

\begin{thm}\label{thm23}\
Let an~ordered Banach algebra $A$ satisfy Axiom~{\bf (A$_1$)}. Let an~element $a\in A_{\rm n}$ have
the~Frobenius normal form. If $r(a)$ is a~pole of $R(\cdot,a)$ then the~residue $a_{-1}$ is
a~finite-rank element.
\end{thm}

{\bf Proof.} The~idea is borrowed from the~proof of the~implication {\bf
(d)}~$\Longrightarrow$~{\bf (c)} of Theorem~2.14 in \cite{Al2}. Let order idempotents
$p_0,p_1,\ldots,p_n$, ${\bf e}=p_n\ge\ldots\ge p_0=0$, determine the~Frobenius normal form of $a$.
Put $q_j=p_jp_{j-1}^{\rm d}$ for $j=1,\ldots,n$. If $r(a)=0$ then $r(a_{q_j})=0$ and $a_{q_j}$ is
irreducible with respect $q_j$ for all $j$, whence $a_{q_j}=0$. Thus, ${\bf
OI}(A)\cap[0,q_j]=\{0,q_j\}$ and so $q_j\gg_{q_j}0$. On the~other hand, $a_{-1}={\bf
e}=\sum\limits_{j=1}^nq_j\in{\cal F}(A)$.

Now we can suppose $r(a)>0$. For the proof, we use induction on $n$. For $n=1$ the~element $a$ is
irreducible and it remains to use Theorem~\ref{thm6}{\bf (c)}. Next, assume that the~desired
assertion is proved if a~parameter of the induction lies between $1$ and $n-1\ge1$. Let us verify
our assertion for $n$.

We consider first the~case of the~identity $r(a_{p_{n-1}})=r(a)$ and show the~inclusion
\begin{equation}\label{23}
(a_{-1})_{p_{n-1}}\in{\cal F}(A).
\end{equation}
The~relation \cite{Al2} $(a_{p_{n-1}})_{-1}=(a_{-1})_{p_{n-1}}$ holds in $A$. If
$r(a)\notin\sigma(a_{p_{n-1}})$ then $(a_{p_{n-1}})_{-1}=0$ and (\ref{23}) is obvious. Let
$r(a)\notin\sigma(a_{p_{n-1}})$. Then \cite{Al2} $r(a)$ is a~pole of
$R_{p_{n-1}}(\cdot,a_{p_{n-1}})$ in $A(p_{n-1})$ (see the~proof of Lemma~\ref{lem21}). The~order
idempotents $p_{n-1},\ldots,p_0$ are \mbox{$a_{p_{n-1}}$-invariant}. If $r(a_{q_j})=r(a_{p_{n-1}})$
then $a_{q_j}$ is irreducible with respect $q_j$ in $A$ and so in $A(p_{n-1})$. By the~induction
hypothesis, the~residue $(a_{p_{n-1}})_{-1}\in{\cal F}(A(p_{n - 1}))$ (see the~proof of
Lemma~\ref{lem21} once more) and, hence, $(a_{p_{n-1}})_{-1}\in{\cal F}(A)$.

Consider the~case of the~identity $r(a_{p_{n-1}^{\rm d}})=r(a)$ and show the~inclusion
\begin{equation}\label{24}
(a_{-1})_{p_{n - 1}^{\rm d}}\in{\cal F}(A).
\end{equation}
Assuming without loss of generality that $r(a)\in\sigma(a_{p_{n-1}}^{\rm d})$, we have
the~relations
$$(a_{-1})_{p_{n - 1}^{\rm d}}=(a_{p_{n - 1}^{\rm d}})_{-1}\in{\cal F}(A(p_{n - 1}^{\rm d}))\subseteq
{\cal F}(A).$$

Both inequalities $r(a_{p_{n-1}})\le r(a)$ and $r(a_{p_{n-1}^{\rm d}})\le r(a)$ cannot be strict
simulta\-neously. To~complete the proof, we consider three possible cases.

{\sf Case~1:} $r(a_{p_{n-1}})=r(a_{p_{n-1}^{\rm d}})=r(a)$. As was shown above, the inclusions
(\ref{23})~and~(\ref{24}) hold. For $\lambda\in\rho(a;A)\setminus\{0\}$, we have
the~identity~\cite{Al2}
\begin{equation}\label{25}
p_{n-1}R(\lambda,a)p_{n-1}^{\rm d}= R(\lambda,a_{p_{n - 1}})p_{n-1}ap_{n-1}^{\rm
d}R(\lambda,a_{p_{n-1}^{\rm d}})
\end{equation}
which implies
$$p_{n - 1}ap_{n-1}^{\rm d}=(a_{p_{n - 1}})_{-1}p_{n-1}ap_{n-1}^{\rm d}(a_{p_{n-1}^{\rm d}})_0+
(a_{p_{n-1}})_0p_{n-1}ap_{n-1}^{\rm d}(a_{p_{n-1}^{\rm d}})_{-1}\in{\cal F}(A)$$ as ${\cal F}(A)$
is an~algebraic ideal. Using the~$a_{-1}$-invariance of $p_{n-1}$, we obtain $a_{-1}\in{\cal
F}(A)$.

{\sf Case~2:} $r(a_{p_{n-1}})=r(a)$ and $r(a_{p_{n-1}^{\rm d}})<r(a)$. Then (\ref{23}) holds.
Moreover, we have~\cite{Al2} $p_{n-1}^{\rm d}a_{-1}=0$. Hence, taking into account~(\ref{25}), we
conclude $a_{-1}\in{\cal F}(A)$.

{\sf Case~3:} $r(a_{p_{n-1}})<r(a)$ and $r(a_{p_{n-1}^{\rm d}})=r(a)$. Then (\ref{24}) holds.
Moreover, we have~\cite{Al2} $a_{-1}p_{n-1}=0$. Hence, taking into account~(\ref{25}), we conclude
$a_{-1}\in{\cal F}(A)$.\hfill$\Box$
\smallskip

We mention the~following important consequence of the~preceding theorem which shows once more that
the~definition of an~$f$-pole is the~right and natural one to use (see~\cite{Al2}, where the~detail
discussion of this notion can be found).

\begin{cor}\label{cor24}\
Let $A$ be a~Dedekind complete and let $a\in A$ be a~spectrally order continuous element with
$r(a)>0$. If $r(a)$ is a~finite-rank pole of $R(\cdot,a)$ then the~residue $a_{-1}$ is
a~finite-rank element.
\end{cor}

{\it It is not known if the~point $r(a)$ is an~$f$-pole of $R(\cdot,a)$ of an~arbitrary irreducible
element $a$ of $A$ such that $r(a)$ is a~pole of $R(\cdot,a)$} (of course, under the~assumptions of
Theorem~\ref{thm6}). As can be shown (see~\cite{Al2}), the~affirmative answer to this question is
equivalent to: {\it $0\le b< a$ implies $r(b)<r(a)$}, which will be discussed in
Section~\ref{sec5}. In~particular, it is not known if the~converse to the~preceding corollary is
true.
\medskip

We now turn our attention to the~conditions of the~primitivity of irreducible element~$a$ in
an~ordered Banach algebra $A$. Recall that an~element $b$ of an~arbitrary Banach algebra $B$ is
called {\it primitive} if the~peripheral spectrum $\sigma_{\rm per}(b)$ contains at most one point;
all other elements of $B$ are called {\it imprimitive}.

We begin with the~next criteria of the~primitivity.

\begin{thm}\label{thm25}\
Let an~ordered Banach algebra $A$ satisfy Axioms {\bf (A$_1$)}-{\bf (A$_5$)} and~{\bf (A$_{\rm
LS}$)}. Let $a\in A$ be a~non-zero irreducible element such that $r(a)$ is a~pole of $R(\cdot,a)$.
The~following statements are equivalent:
\begin{description}
\item[(a)] The~element $a$ is primitive;
\item[(b)] The~element $a^m$ is primitive for all $m\in{\Bbb N}$;
\item[(c)] The~element $a^m$ is irreducible for all $m\in{\Bbb N}$;
\item[(d)] The~sequence $\{(\frac{a}{r(a)})^n\}$ converges to an~algebraically strictly positive
element.
\end{description}
\end{thm}

{\bf Proof.} The~implication {\bf (a)}~$\Longrightarrow$~{\bf (b)} follows at once from
the~identity $\sigma(a^m)=f(\sigma(a))$ with $f(z)=z^m$ and the~implication {\bf
(b)}~$\Longrightarrow$~{\bf (a)} is obvious.

{\bf (a)}~$\Longrightarrow$~{\bf (c)} If $a^m$ is reducible for some $m\in{\Bbb N}$ then $m>1$ and,
in view of Theorem~\ref{thm12}, the~points $r(a)e^{i\frac{2\pi}{m''}j}\in\sigma(a)$ for all
$j=0,1,\ldots,m''-1$ with $m''>1$, a~contradiction.

{\bf (c)}~$\Longrightarrow$~{\bf (a)} Proceeding by contradiction and using Axiom~{\bf (A$_{\rm
LS}$)} and Theorem~\ref{thm20}, we conclude the~validity of the~identity $\sigma_{\rm
per}(a)=r(a)H_k$ with $k>1$. Taking into account Theorem~\ref{thm12} once more, we obtain
the~reducibility $a^m$ of for some $m>1$, which is impossible.

{\bf (a)}~$\Longrightarrow$~{\bf (d)} We shall prove first the~next assertion: {\it If $b$ is
a~primitive positive element of an~arbitrary ordered Banach algebra $A_0$ such that $r(b)>0$ is
a~simple pole of $R(\cdot,b)$ then $(\frac{b}{r(b)})^n$ converges to the~residue $b_{-1}$ of
$R(\cdot,b)$ at $r(b)$}. Indeed, we define the~element~$c$ by $c=b-r(b)b_{-1}$. The~Spectral
Mapping Theorem yields $\sigma(c)=(\sigma(c)\cup\{0\})\setminus\{r(b)\}$. Taking into account
the~primitivity of $b$, we infer $r(c)<r(b)$ or $r(\frac{c}{r(b)})<1$. By the~Gelfand formula,
the~equality $r(\frac{c}{r(b)})=\lim\limits_{n\to\infty}\frac{\|c^n\|_{A_0}}{r(b)^n}$ holds and,
in~particular, $(\frac{c}{r(b)})^n\to0$ as $n\to\infty$. Using the~identities
$bb_{-1}=b_{-1}b=r(b)b_{-1}$, we obtain the~relations $b_{-1}c=cb_{-1}=0$. Therefore,
$b^n=c^n+r(b)^nb_{-1}$ for all $n$. Consequently, the~relations
$(\frac{b}{r(b)})^n=(\frac{c}{r(b)})^n+b_{-1}\to b_{-1}$ hold, as~required. Now it only remains to
remember Theorem~\ref{thm6}{\bf (b)},{\bf (c)}. According to it, $r(a)$ is a~simple pole of
$R(\cdot,A)$ and $a_{-1}\gg0$ as $a$ is irreducible.

{\bf (d)}~$\Longrightarrow$~{\bf (c)} If $a^m$ is reducible for some number $m>1$ then $p^{\rm
d}a^mp=0$ for some non-trivial $p\in{\bf OI}(A)$. Evidently, $p^{\rm d}a^{mn}p=0$ for all
$n\in{\Bbb N}$. Thus, $0=p^{\rm d}(\frac{a}{r(a)})^{mn}p\to p^{\rm d}a_{-1}p$ as $n\to\infty$ and,
hence, $p^{\rm d}a_{-1}p=0$. The~latter contradicts to the~algebraic strict positivity of
$a_{-1}$.\hfill$\Box$
\smallskip

A~non-zero order idempotent $p$ of $A$ is called {\it order minimal} if the~equality
$\sum\limits_{j=1}^np_j={\bf e}$, where $p_j\in{\bf OI}(A)$ for $j=1,\ldots,n$ and
$p_{j'}p_{j''}=0$ for $j'\neq j''$, implies the~existence of an~index $j_0$ satisfying $p\le
p_{j_0}$.

\begin{cor}\
Under the~assumptions of Theorem~{\rm \ref{thm25}}, each of the~following conditions guarantees
the~primitivity of the~element $a$:
\begin{description}
\item[(a)] The~element $a\gg0$;
\item[(b)] The~element $a_p>0$ for some order minimal $p\in{\bf OI}(A)$.
\end{description}
\end{cor}

{\bf Proof.} {\bf (a)} We shall show the~irreducibility of $a^m$ for all $m\in{\Bbb N}$. In view of
part~{\bf (a)} of the~preceding theorem, the~latter implies the~desired assertion. To this end, we
assume $q^{\rm d}a^mq=0$ for some $m>1$ and $q\in{\bf OI}(A)$. Since $A$ has a~disjunctive product,
there exists an~order idempotent $q_1$ satisfying $q^{\rm d}aq_1=q_1^{\rm d}a^{m-1}q=0$. We can
suppose $q\neq{\bf e}$. Then $q_1=0$ and so $a^{m-1}q=0$. If $m>2$ then there exists an~order
idempotent $q_2$ satisfying $aq_2=q_2^{\rm d}a^{m-2}q=0$ and so $a^{m-2}q=0$. Finally, using
a~reverse finite induction, we obtain $aq=0$ or $q=0$, as required.

{\bf (b)} Proceeding by contradiction and taking into account Theorems \ref{thm25} and~\ref{thm12},
we find $p_1,\ldots,p_m\in{\bf OI}(A)$ with $m>1$ determining the~cyclic form of the~element~$a$
and, in~particular, $p_ja=ap_{j+1}$ for all $j=1,\ldots,m$. In view of the~order minimality of $p$,
we choose an~index~$j_0$ satisfying $p\le p_{j_0}$. Then the~equalities $pa=pp_{j_0}a=pap_{j_0+1}$
hold and, hence, $a_p=pap_{j_0+1}p=0$, a~contradiction.\hfill$\Box$
\smallskip

If $b$ is an~element of a~Banach algebra $B$ with a~unit ${\bf u}$ such that $r(b)\in\sigma(b)$
then $b+\lambda{\bf u}$ is primitive for all numbers $\lambda>0$. The~next result makes more
precisely this fact.

\begin{cor}\
Let $A$ satisfy Axioms {\bf (A$_1$)}-{\bf (A$_5$)} and let $a,b\in A^+$. If $a+b$ is irreducible,
$b_p>0$ for all $0<p\in{\bf OI}(A)$, and $r(a+b)$ is a~pole of $R(\cdot,a+b)$ then $a+b$ is
primitive.
\end{cor}

{\bf Proof.} Again, proceeding by contradiction and taking into account Theorems \ref{thm25}
and~\ref{thm12}, we find $p_1,\ldots,p_m\in{\bf OI}(A)$ with $m>1$ determining the~cyclic form
of~$a+b$ and, in~particular, $p_j(a+b)=(a+b)p_{j+1}$ for all $j=1,\ldots,m$. Therefore,
$0=p_j(a+b)p_j=a_{p_j}+b_{p_j}>0$, a~contradiction.\hfill$\Box$
\smallskip

An~element $a$ of $A$ is said to be {\it symmetric} whenever $paq=qap$ for all $p,q\in{\bf OI}(A)$.

\begin{cor}\
Suppose that all assumptions of Theorem~{\rm \ref{thm25}} are satisfied and, in~addi\-tional,
the~element $a$ is symmetric. Then the~element $a$ is primitive if and only if $a+r(a){\bf e}$ is
invertible.
\end{cor}

{\bf Proof.} The~necessity is obvious. We shall prove the~sufficiency. Proceeding by
contra\-diction, we find $p_1,\ldots,p_m\in{\bf OI}(A)$ with $m>1$ determining the~cyclic form
of~$a$. In~view of the~validity of the~implication {\bf (b)}~$\Longrightarrow$~{\bf (c)} of
Theorem~\ref{thm12} and our condition, the~natural number $m$ is odd. Then $p_1ap_2=ap_2>0$. Since
$a$ is symmetric, $p_2ap_1>0$. On the~other hand, $p_2ap_1=ap_3p_1=0$, a~contradiction.\hfill$\Box$
\smallskip

As~follows from Theorem~\ref{thm25}, if an~element $a$ of $A$ is irreducible and primitive then
$a^m$ is irreducible for all $m\in{\Bbb N}$. Unfortunately, the~relation $a^m\gg0$ need not hold
for any $m$ in the~case of an~integral operator even (see~\cite{Al}). Nevertheless, we have
the~next result.

\begin{cor}\
Suppose that all assumptions of Theorem~{\rm \ref{thm25}} are satisfied and, in~addi\-tional,
the~element $a$ is primitive. If $a^{m_0}\gg0$ for some $m_0\in{\Bbb N}$ then $a^m\gg0$ for all
natural $m\ge m_0$.
\end{cor}

{\bf Proof.} Proceeding by contra\-diction, we find $m>m_0$ and non-zero $p,q\in{\bf OI}(A)$
satisfying $pa^mq=0$. Obviously, $pa^{m_0}a^{m-m_0}q=0$. Since the~algebra $A$ has a~disjunctive
product, we have $pa^{m_0}p_1=p_1^{\rm d}a^{m-m_0}q=0$ for some $p_1\in{\bf OI}(A)$. According to
our condition, $p_1=0$ and so $a^{m-m_0}q=0$. In view of Theorem~\ref{thm25}, the~element
$a^{m-m_0}$ is irreducible and, hence, $q=0$, a~contradiction.\hfill$\Box$
\smallskip

For the~case of an~arbitrary (not necessarily irreducible) element $a\in A^+$, we have the~next.

\begin{prop}\
Let $A$ be an~ordered Banach algebra. Let $a\in A$ be a~positive element such that $r(a)$ is a~pole
of $R(\cdot,A)$ and every point $\alpha\in\sigma_{\rm per}(a)$ is an~eigenvalue of either
the~operator $L_a$ or the~operator $R_a$ on $A$. The~following statements hold:
\begin{description}
\item[(a)] The~sequence $a^n\to0$ as $n\to\infty$ if and only if $r(a)<1$;
\item[(b)] If $r(a)=1$ then the~sequence $\{a^n\}$ is convergent if and only if $r(a)$ is a~simple
pole of $R(\cdot,a)$ and the~element $a$ is primitive;
\item[(c)] If $r(a)>1$ then the~sequence $\{a^n\}$ is not convergent.
\end{description}
\end{prop}

{\bf Proof.} {\bf (b)} Let $r(a)=1$ and let $\{a^n\}$ be convergent. Then $\|a^n\|_A\le c$ for all
$n$ and some constant $c\in{\Bbb R}^+$. For all $\lambda\in{\Bbb R}$, $\lambda>1$, we have
$$(\lambda-1)^2\|R(\lambda,a)\|_A=
(\lambda-1)^2\Big\|\sum\limits_{j=0}^\infty\frac1{\lambda^{j+1}}a^j\Big\|_A\le$$
$$(\lambda-1)^2\sum\limits_{j=0}^\infty\frac1{\lambda^{j+1}}\|a^j\|_A\le
C(\lambda-1)^2\sum\limits_{j=0}^\infty\frac1{\lambda^{j+1}}=
C(\lambda-1)^2\frac\lambda{\lambda-1}=C\lambda(\lambda-1)\to0$$ as $\lambda\downarrow1$. Thus,
$r(a)$ is a~simple pole of $R(\cdot,a)$. Now we consider $\alpha\in\sigma_{\rm per}(a)$. In view of
our condition, $\alpha$ is an~eigenvalue of $L_a$ (the~case of $R_a$ is analogous), i.e.,
$ab=\alpha b$ for some non-zero $b\in A$. Obviously, $a^nb=\alpha^nb$ for all $n\in{\Bbb N}$ and
the~sequence $\{\alpha^nb\}$ converges. Therefore, $\{\alpha^n\}$ converges. The~latter is possible
for the~case of $\alpha=1$ only and, hence, $a$ is primitive. The~converse assertion, namely,
the~relation $a^n\to a_{-1}$ as $n\to\infty$, was shown in the~proof of the~implication {\bf
(a)}~$\Longrightarrow$~{\bf (d)} of Theorem~\ref{thm25}.

{\bf (a)} The~sufficiency is clear. We shall check the~necessity. Assume that $a^n\to0$ as
$n\to\infty$. This implies $r(a)\le1$. If $r(a)=1$ then, as was mentioned above, $a^n\to a_{-1}$
and, hence, $a_{-1}=0$, a~contradiction.

{\bf (c)} If $r(a)=1$ then $\{a^n\}$ is unbounded and, in~particular, is not
convergent.\hfill$\Box$

\section{Other viewpoints on the~irreducibility and\newline the~primitivity}

Recall that an~element $a\in A^+$ is said to be irreducible whenever the~equality $p^{\rm d}ap=0$,
where $p\in{\bf OI}(A)$, implies $p=0$ or $p={\bf e}$. Irreducible elements were introduced
in~\cite{BGr} for the~case of Banach lattice algebras and in~\cite{Al2} for the~case of ordered
Banach algebras. Under the~natural assumptions, these elements have nice spectral properties and,
under such a~notion of irreducibility, the~theorem about the~Frobenius normal form holds.
Nevertheless, the~purpose of this section is to establish some results which allow us in a~new
fashion to take a~glance at the~algebraic nature of such notions as the~irreducibility and
the~primitivity and, thus, to break some ground for further research of these notions in an~ordered
Banach algebras.
\medskip

Let $E$ be a~Dedekind complete Banach lattice. As usual, we denote the~set of all opera\-tors on
$E$ of the~form $\sum\limits_{j=1}^nf_j\otimes x_j$, where $f_j\in E_{\rm n}^\sim$, $x_j\in E$, and
$(f_j\otimes x_j)x=(f_j(x))x_j$ for all $j=1,\ldots,n$ and $x\in E$, by $E_{\rm n}^\sim\otimes E$.
The~band $(E_{\rm n}^\sim\otimes E)^{\rm dd}$ generated by $E_{\rm n}^\sim\otimes E$ in the~Banach
lattice $L_{\rm r}(E)$ with the~$r$-norm $\|T\|_r=\|T\|_{B(E)}$ is called \mbox{\cite[p.
193]{AbrAl}} the~{\it band of abstract integral operators}. As is easy to see, $(E_{\rm
n}^\sim\otimes E)^{\rm dd}\subseteq L_{\rm n}(E)$ and if $T\in L_{\rm r}(E)$ and $S\in(E_{\rm
n}^\sim\otimes E)^{\rm dd}$ then $TS,ST\in(E_{\rm n}^\sim\otimes E)^{\rm dd}$. In~particular, under
the~$r$-norm, $(E_{\rm n}^\sim\otimes E)^{\rm dd}$ is an~ordered Banach algebra (possibly, without
a~unit). If $E$~is a~function space and $E_{\rm n}^\sim$ separates points of $E$ then, by
the~Lozanovsky theorem \mbox{\cite[p. 199]{AbrAl}}, $(E_{\rm n}^\sim\otimes E)^{\rm dd}$ coincides
with the~band of regular integral operators on $E$. Next, as can be shown (see~\cite{Sch2}), if
an~arbitrary Banach lattice $E$ possesses the~non-trivial band $E_{\rm c}^\sim$ of all
$\sigma$-order continuous functionals, i.e., $E_{\rm c}^\sim\neq\{\emptyset\}$, and admitts
a~$\sigma$-order continuous irreducible operator $T$ then $E_{\rm n}^\sim=E_{\rm c}^\sim$ and
$E_{\rm n}^\sim$ separates the~points of $E$. We also recall that in a~Banach algebra $B$ with or
without a~unit the~{\it wedge operator}~$W_b$ on $B$ \mbox{\cite[pp. 17, 70]{BMSmW}}, where $b\in
B$, is defined by $W_b=bab$ for $a\in B$. As is easy to see, $r(W_b)\le r(b)^2$ (we put
$r(b)=\lim\limits_{n\to\infty}\|b^n\|_B^{\frac1n}$ if $B$ does not have a~unit).

\begin{thm}\label{thm31}\
Let $E$ be a~Dedekind complete Banach lattice such that $E_{\rm n}^\sim$ separates points of $E$.
Let $T$ be a~non-zero positive order continuous operator on $E$ such that $r(T)$ is a~pole of
$R(\cdot,T)$. If the~restriction of the~wedge operator $W_T$ to $(E_{\rm n}^\sim\otimes E)^{\rm
dd}$ is an~irreducible operator then two operators $T$ and $T'$ are also irreducible and, moreover,
primitive, where $T'$ is the~restriction of the~adjoint operator $T^*$ to $E_{\rm n}^\sim$.
\end{thm}

{\bf Proof.} Consider \mbox{a~$T$-invariant} band $B\neq\{0\}$. Fix $\lambda>r(T)$. Using
\mbox{the~$R(\lambda,T)$-invariance} of this band and the~inequality $TR(\lambda,T)\le\lambda
R(\lambda,T)$, we find a~non-zero element $z_0\in B^+$ satisfying $Tz_0\le\lambda z_0$. For
an~arbitrary non-zero functional $h_0\in(E^*)^+$ such that $T^*h_0\le\lambda h_0$, we have
$W_T(h_0\otimes z_0)\le\lambda^2h_0\otimes z_0$. As is easy to see, the~operator $W_T$ is order
continuous on~$(E_{\rm n}^\sim\otimes E)^{\rm dd}$. Therefore, $h_0\otimes z_0$ is a~weak order
unit in $(E_{\rm n}^\sim\otimes E)^{\rm dd}$. If $z\in E^+$ and $z_0\perp z$ then $(h\otimes
z_0)\wedge(h\otimes z)=h\otimes(z_0\wedge z)=0$, whence $h\otimes z=0$ or $z=0$. Thus, $z_0$ is
a~weak order unit. Consequently, $B=E$ and the~irreducibility of $T$ has been proved.
In~particular, $T'$~is (see~\cite{Al}, the~proof of Theorem~1) also irreducible. Moreover, there
exist a~weak order unit $x_0\in E$ and a~strictly positive functional $f\in E_{\rm n}^\sim$
satisfying $Tx_0=r(T)x_0$ and $T^*f_0=r(T)f_0$.

Now let us verify the~primitivity of the~operator $T$. In view of the~inclusion \mbox{\cite[p.
256]{AbrAl}} $\rho_\infty(T^*;B(E^*))\subseteq\rho(T';B(E_{\rm n}^\sim))$, the~latter implies
the~primitivity of $T'$. Proceeding by contradiction, we find (\cite{Al}; see also part~{\bf (b)}
of Theorem~\ref{thm12}) elements $y_1,\ldots,y_m$, where $m>1$, satisfying
$\sum\limits_{j=1}^my_j=x_0$, $y_{j'}\wedge y_{j''}=0$ for $j'\neq j''$, and $Ty_{j+1}=r(T)y_j$ for
$j=1,\ldots,m$. Define the~functionals $g_j\in E^*$ via the~formula $g_j=P_{y_j}^*f_0$, where
$P_{y_j}$ is the~order projection from $E$ onto the~projection band $B_{y_j}$ generated by $y_j$.
Obviously, $\sum\limits_{j=1}^mg_j=f_0$, $g_{j'}\perp g_{j''}$ for $j'\neq j''$, and $g_j\in E_{\rm
n}^\sim$. Taking into account the~identities $T^*g_j=T^*P_{y_j}^*f_0=P_{y_{j+1}}^*T^*f_0$, we have
the~relation $T^*g_{j'}\perp T^*g_{j''}$ for $j'\neq j''$. Moreover, $(T^*g_j)y_k=g_j(y_{k-1})=0$
for $k\neq j+1$, whence $(T^*g_j)(x_0-y_{j+1})=g_ky_{j+1}=0$ and so $(T^*g_j\wedge g_k)x_0=0$.
Since $x_0$ is a~weak order unit, the~last equality implies $T^*g_j\perp g_k$ for $k\neq j+1$.
Consequently,
$$T^*g_j-r(T)g_{j+1}\perp r(T)\sum\limits_{\stackrel{n=1}{n\neq j+1}}^mg_n-
\sum\limits_{\stackrel{n=1}{n\neq j}}^mT^*g_n.$$ On the~other hand,
$\sum\limits_{n=1}^mT^*g_n=r(T)\sum\limits_{n=1}^mg_n$ or
$$T^*g_j-r(T)g_{j+1}=r(T)\sum\limits_{\stackrel{n=1}{n\neq j+1}}^mg_n-
\sum\limits_{\stackrel{n=1}{n\neq j}}^mT^*g_n,$$ whence $T^*g_j=r(T)g_{j+1}$ for all
$j=1,\ldots,m$. Put $f_j=g_{m-j+1}$. Evidently,
$$T^*f_{j+1}=T^*g_{m-j}=r(T)g_{m-j+1}=r(T)f_j.$$
Thus, $W_T\sum\limits_{n=1}^mf_n\otimes y_n=r(T)^2\sum\limits_{n=1}^mf_n\otimes y_n$. Therefore,
$\sum\limits_{n=1}^mf_n\otimes y_n$ is a~weak order unit in the~band~$(E_{\rm n}^\sim\otimes
E)^{\rm dd}$. This contradicts to the~relation $f_j\otimes y_k\perp\sum\limits_{n=1}^mf_n\otimes
y_n$ for all $j\neq k$.\hfill$\Box$
\smallskip

It is not known if the~converse to the~assertion of the~preceding theorem holds. Nevertheless, as
the~next result shows, we can assert the~converse in the~case of the~Banach algebra $M_n({\Bbb C})$
of all $n\times n$ matrices with complex entries and the~natural multiplication and order. We
mention first that a~matrix $A\in M_n({\Bbb C})$ is irreducible if and only if the~transpose $A^t$
of $A$ is irreducible.

\begin{thm}\
Let $T$ be an~$n\times n$~positive matrix. Then $T$ is irreducible and primitive if and only if
the~wedge operator $W_T$ on $M_n({\Bbb C})$ is irreducible.
\end{thm}

{\bf Proof.} In view of Theorem~\ref{thm31}, it is enough to verify the~necessity. Clearly, we can
assume~$n>1$. Let the~inequality $W_TS\le\lambda S$ hold for a~number $\lambda\ge0$ and a~non-zero
positive matrix $S\in M_n({\Bbb C})$. Since the~matrix $T$ is primitive, $T^k$ is strongly positive
for some $k\in{\Bbb N}$ (see Theorem~\ref{thm25}), i.e., all its entries are strictly positive.
Therefore, $\lambda^kS\ge W_T^kS=T^kST^k$ and so $S$ is also strongly positive. Finally, $W_T$ is
irreducible.\hfill$\Box$
\medskip

The~following result suggests another approach to the~notion of irreducibility in ordered Banach
algebras.

\begin{thm}\
Let $T$ be an~$n\times n$~matrix. Then $T$ is irreducible if and only if the~operator $L_T+R_T$ on
$M_n({\Bbb C})$ is irreducible.
\end{thm}

{\bf Proof.} Necessity. We recall first that (see the~equalities~(\ref{7})) $(L_T+R_T)Q=TQ+QT$,
where $Q\in M_n({\Bbb C})$. Let the~inequality $L_TS\le\lambda S$ holds for a~number $\lambda\ge0$
and a~non-zero positive matrix $S\in M_n({\Bbb C})$, $S=[s_{ij}]$. There exist indexes
$i_0,j_0=1,\ldots,n$ such that the~entry $s_{i_0j_0}>0$. If $s^{j_0}$ is the~$j_0^{\rm th}$ column
of $S$ then, in view of the~inequality $TS\le\lambda S$, we have $Ts^{j_0}\le\lambda s^{j_0}$.
Taking into account the~irreducibility of $T$, we obtain $s_{ij_0}>0$ for all $i=1,\ldots,n$. On
the~other hand, the~inequality $S^tT^t\le\lambda S^t$ holds, $T^t$~is irreducible, and
$s_{j_0i}^{t_0}>0$ for all $i$, where $S^t=[s_{ij}^t]$. As was shown above, $s_{ji}^t>0$ or
$s_{ij}>0$ for all $i,j$, and we are done.

Sufficiency. Let the~operator $L_T+R_T$ be irreducible. In~particular, $L_T+R_T\ge0$. Then
$0\le(L_T+R_T)I=2T$ or $T\ge0$. Assume that $T$ is reducible. Then for some $k=1,\ldots,n-1$ there
exist indexes $j_1,\ldots,j_k$ satisfying $t_{ij}=0$ for all $i\notin J$ and $j\in J$, where
$J=\{j_1,\ldots,j_k\}$. Consider the~band $B$ in $M_n({\Bbb C})$ defined by
$$B=\{S\in M_n({\Bbb C}):s_{ij}=0 \ \text{for all} \ i\notin J \ \text{and all} \ j\}.$$
Let $S\in B$. If $TS=[(ts)_{ij}]$ then for $i\notin J$, we have
$$(ts)_{ij}=\sum\limits_{m=1}^nt_{im}s_{mj}=
\sum\limits_{m\notin J}t_{im}s_{mj}+\sum\limits_{m\in J}t_{im}s_{mj}=0.$$ If $ST=[(st)_{ij}]$ then
$(st)_{ij}=0$ for $i\notin J$. Thus, the~band $B$ is \mbox{$L_T+R_T$-invariant}, which is
impossible.\hfill$\Box$

\section{When does $0\le b< a$ imply $r(b)<r(a)$?}\label{sec5}

Let $T$ be an~order continuous irreducible operator on a~Banach lattice $E$, let $r(T)$~be a~pole
of $R(\cdot,T)$ of order $k$, and let the~coefficient $T_{-k}$ of the~Laurent series expansion
of~$R(\cdot,T)$ around~$r(T)$ also be order continuous (the~latter holds if, e.g., the~Lorenz
seminorm~$\|\cdot\|_L$ on~$E$ is~a~norm). Then \cite{Al2} the~operator inequalities $0\le S< T$,
where $S,T\in B(E)$, imply the~spectral radius inequality $r(S)<r(T)$. In~particular, if $S,T\in
M_n({\Bbb C})$, $T$ is irreducible, and $0\le S< T$ then $r(S)<r(T)$. As was mentioned above (see
remarks after Corollary~\ref{cor24}), the~analogous question, i.e., the~validity of the~inequality
$r(b)<r(a)$ where $0\le b< a$ and the~element $a$ is irreducible, remains open in the~case of
an~ordered Banach algebra. The~purpose of the~present section is to discuss a~number of additional
conditions under which the~answer to this question is affirmative. We mention at once that in
research of this problem the~assumption that $r(a)$ is a~pole of $R(\cdot,a)$ is natural absolutely
and cannot be even reject in the~case of operators (see~\mbox{\cite{Al, Al4}}).

\begin{thm}\label{thm34}\
Let $A$ be a~finite-dimensional ordered Banach algebra with a~disjunctive product and let $a,b\in
A$. If the~element $a$ is irreducible then $0\le b<a$ implies $r(b)<r(a)$.
\end{thm}

{\bf Proof.} We claim that all assumptions of Theorem~\ref{thm6} hold. Indeed, if elements
$q_j\in{\bf OI}(A)$ with $j=1,\ldots,k$ and $q_{j'}q_{j''}=0$ for $j'\neq j''$ then
$q_1,\ldots,q_k$ are linearly independent. Therefore, there exists a~maximal collection of pairwise
disjoint elements $\{p_1,\ldots,p_n\}$ in ${\bf OI}(A)$. Obviously, $\sum\limits_{j=1}^np_j={\bf
e}$. We have the~identity
\begin{equation}\label{26}
{\bf OI}(A)=\Big\{\sum\limits_{j\in J}p_j:J\subseteq\{1,\ldots,n\}\Big\}.
\end{equation}
Actually, if an~order idempotent $p$ of $A$ satisfies $0<pp_j<p_j$ for some $j=1,\ldots,n$ then
$0<pp_j<p_j$ and the~collection $\{p_1,\ldots,p_{j-1},pp_j,p^{\rm d}p_j,p_{j+1},\ldots,p_n\}$
consists of pairwise disjoint elements. The~latter contradicts to the~maximality of
$\{p_1,\ldots,p_n\}$ and, hence, either $pp_j=0$ or $pp_j=p_j$. Putting $J=\{j:pp_j>0\}$, we obtain
$p=p\sum\limits_{j=1}^np_j=\sum\limits_{j\in J}pp_j=\sum\limits_{j\in J}p_j$, and (\ref{26})~has
been checked. In~particular, the~Boolean algebra ${\bf OI}(A)$ is Dedekind complete. Consider a~net
$\{b_\alpha\}$ satisfying $b_\alpha\downarrow0$ in ${\bf OI}(A)$. As~was shown above,
$b_\alpha=\sum\limits_{j\in J_\alpha}p_j$ for every $\alpha$, where
$J_\alpha\subseteq\{1,\ldots,n\}$. Evidently, if $b_{\alpha'}<b_{\alpha''}$ then
$J_{\alpha'}\subsetneq J_{\alpha''}$. From the~latter, we conclude easily the~existence of
an~index~$\alpha_0$ such that $b_{\alpha_0}=0$. Thus, the~equality $A^+=A_{\rm n}$ holds. On
the~other hand, as is well known, the~spectrum $\sigma(x;B)$ of every element $x$ of an~arbitrary
finite-dimensional Banach algebra $B$ with a~unit is finite and consists of poles of
the~resolvent~$R(\cdot,x)$. Now it only remains to use part~{\bf (d)} of
Theorem~\ref{thm6}.\hfill$\Box$

As the~next example shows, the~preceding theorem is not valid without the~assump\-tion about
a~disjunctive product.

\begin{exm}\
{\rm Consider the~ordered Banach algebra $A_0={\Bbb C}^n$, where $n\in{\Bbb N}$ and $n\ge2$, under
the~natural algebraic operations, multiplication, and order and under some algebra-norm. Then
the~ordered Banach algebra $A$ obtained from $A_0$ by adjoining a~unit gives the~required example.
Another example is the~ordered Banach algebra $A={\Bbb C}^2$ under the~natural algebraic operations
and order, the~multiplication given by $(x_1,x_2)(y_1,y_2)=(x_1y_1,x_1y_2+x_2y_1)$, and the~norm
$\|(x_1,x_2)\|_A=|x_1|+|x_2|$. Indeed, as is easy to see, the~element ${\bf e}=(1,0)$ is a~unit of
$A$, ${\bf OI}(A)=\{0,{\bf e}\}$, and the~spectrum $\sigma(x;A)=\{x_1\}$, where $x=(x_1,x_2)$. Then
every element of the~algebra $A$ is irreducible, $0\le(1,0)<(1,1)$, and $r((1,0))=r((1,1))=1$.

The~situation does not change in the~case of the~algebra of the~form $B(E)$. Indeed, let $H$~be
an~arbitrary real Hilbert space and let $z\in H$ with $\|z\|_H=1$. Under the~order generated by
the~ice cream cone $K=\{y\in H:\langle y,z\rangle\ge\frac1{\sqrt{2}}\|y\|_H\}$, the~space $H$ is
a~real ordered Banach space. Since the~cone $K$ is generating, $B(H)$ is a~real ordered Banach
algebra and, hence, $B(H_{\Bbb C})$ is a~complex ordered Banach algebra, where $H_{\Bbb C}$ is
the~complexification of $H$. If $\dim{H}\ge3$ then \cite{Al3} the~center $(B(H_{\Bbb
C}))_I=\{\lambda I:\lambda\in{\Bbb R}\}$ and, in~particular, we have ${\bf OI}(B(H_{\Bbb
C}))=\{0,I\}$. There exists (see~\cite{Al3} once more) a~non-zero positive operator $T$ satisfying
$T^2=0$. Obviously, $I+T$ is an~irreducible element of $B(H_{\Bbb C})$, $0\le I<I+T$, and
$r(I)=r(I+T)=1$. On the~other hand, if $E$ is a~two-dimensional ordered Banach space with
generating cone $E^+$ then, as is well known, $E^+$~is a~lattice cone. From the~latter follows
easily that the~ordered Banach algebra $B(E)$ has a~disjunctive product and, hence,
Theorem~\ref{thm34} can be applied in this case.}\hfill$\Box$
\end{exm}

An~arbitrary Banach algebra $B$ with a~unit is said to be {\it Fredholm} if the~following three
conditions hold:
\begin{description}
\item[(a)] The~open subset $\Phi(B)$ of $B$ is determined and $\Phi(B)=-\Phi(B)$.
Elements of~$\Phi(B)$ are called {\it Fredholm};
\item[(b)] Two functions ${\rm nul}, {\rm def}:B\to{\Bbb N}\cup\{0,\pm\infty\}$ are determined
such that ${\rm nul}\! \ b={\rm nul}\! \ (-b)$ and ${\rm def}\! \ b={\rm def}\! \ (-b)$ for all
$b\in B$ and the~set ${\rm Inv}\! \ B$ of invertible elements of $B$ satisfies ${\rm Inv}\! \
B=\{b\in B:{\rm nul}\! \ b={\rm def}\! \ b=0\}$;
\item[(c)] The~{\it punctured neighbourhood property} holds: if $b\in\Phi(B)$ then there exists
a~number $\epsilon>0$ such that ${\rm nul}\! \ (\lambda-b)$ and ${\rm def}\! \ (\lambda-b)$ are
constants on the~set $\{\lambda\in{\Bbb C}:0<|\lambda|<\epsilon\}$.
\end{description}

As is well known, the~Banach algebra $B(Z)$, where $Z$ is an~arbitrary Banach space, is Fredholm
(see~\mbox{\cite[Section~4.4]{AbrAl}}). In this case, $\Phi(B(Z))$ coincides with the~class of all
Fredholm operators on $Z$ and ${\rm nul}\! \ T=\dim{N(T)}$ and ${\rm def}\! \ T={\rm codim}\! \
R(T)$ for all $T\in\Phi(B(Z))$. Moreover, every Banach algebra $B$ with a~unit ${\bf u}$ is
Fredholm (see~\mbox{\cite[Sections F2 anf F3]{BMSmW}}), i.e., the~set $\Phi(B)$ and two functions
${\rm nul}$ and ${\rm def}$ on $B$ satisfying the~required properties can be defined. Nevertheless,
in~\cite{BMSmW} for the~deter\-mi\-nation of these objects the~notion of {\it inessential ideal}
$J$ was used. That is, $J$ is an~algebraic ideal and zero is the~only possible accumulation point
of $\sigma(b;B)$ for each $b\in J$. In this case, the~set $\Phi(B)$ is defined by
\begin{equation}\label{27}
\Phi(B)=\{b\in B:{\bf u}-ab,{\bf u}-ba\in J \ \text{for some} \ a\in B\}.
\end{equation}
In the~definition of the~Fredholm algebra given above the~notion of inessential ideal is not
required. Moreover, in some cases, e.g., of the~algebra $C(K)$, where the~compact space $K$ does
not contain unisolated points, the~zero ideal is a~unique inessential ideal.

\begin{thm}\
Let an~ordered Banach algebra $A$ satisfy Axiom~{\bf (A$_1$)} and let the~cone $A^+$ be normal. Let
$a,b\in A$ be such that $0\le b<a$, the~element $a$ is irreducible, $r(a)$~is a~pole of
$R(\cdot,a)$. Then each of the~following conditions guarantees the~inequality $r(b)<r(a)$:
\begin{description}
\item[(a)] $A$ is a~Fredholm algebra and $r(a)-a\in\Phi(A)$;
\item[(b)] There exists an~inessential ideal $J$ of $A$ which contains the~residue $a_{-1}$;
\item[(c)] The~ideal ${\cal F}(A)$ of finite-rank elements is inessential.
\end{description}
\end{thm}

{\bf Proof.} {\bf (a)} Proceeding by contradiction, we assume $r(b)=r(a)$. For arbitrary
$\epsilon\in[0,1)$, we define the~element $a_\epsilon$ by $a_\epsilon=(1-\epsilon)a+\epsilon b$.
Obviously, $0\le b\le a_\epsilon\le a$. Therefore, in~view of the~normality of the~cone $A^+$,
$r(a_\epsilon)=r(a)$. Moreover, $r(a)-a_\epsilon\to r(a)-a$ as $\epsilon\to0$ and, hence,
$r(a)-a_\epsilon\in\Phi(A)$ for sufficiently small $\epsilon$. Fix such a~number $\epsilon$. Then
${\rm nul}\! \ (\lambda-a_\epsilon)$ and ${\rm def}\! \ (\lambda-a_\epsilon)$ are constants on
the~set $\{\lambda\in{\Bbb C}:0<|\lambda-r(a)|<\delta\}$ for some $\delta>0$. On the~other hand,
for $\lambda$ close to $r(a)$, we have the~inclusion $\lambda-a\in{\rm Inv}\! \ A$ or ${\rm nul}\!
\ (\lambda-a)={\rm def}\! \ (\lambda-a)=0$. Since the~element $\lambda-a_\epsilon$ is invertible
for $|\lambda|>r(a)$, the~identities ${\rm nul}\! \ (\lambda-a_\epsilon)={\rm def}\! \
(\lambda-a_\epsilon)=0$ hold for $|\lambda|>r(a)$ and so for all numbers $\lambda$ satisfying
$0<|\lambda-r(a)|<\delta$. Therefore, for such $\lambda$ the~element $\lambda-a_\epsilon$ is
invertible and, in~particular, $r(a)$ is an~isolated point of $\sigma(a_\epsilon)$. Using
the~normality cone $A^+$ once more, we conclude that $r(a)$ is a~simple pole of
$R(\cdot,a_\epsilon)$. The~residue $(a_\epsilon)_{-1}$ satisfies the~relations
$$
 r(a)(a_\epsilon)_{-1}=a_\epsilon(a_\epsilon)_{-1}\le a(a_\epsilon)_{-1} \ \ \text{or} \ \
 0\le(a-r(a))(a_\epsilon)_{-1}.
$$
On the~other hand, $a_{-1}(a-r(a))(a_\epsilon)_{-1}=0$. Taking into account the~relation
$a_{-1}\gg0$, we get
$$a(a_\epsilon)_{-1}=r(a)(a_\epsilon)_{-1}=a_\epsilon(a_\epsilon)_{-1}=
((1-\epsilon)a+\epsilon b)(a_\epsilon)_{-1}$$ and, hence, $(a-b)(a_\epsilon)_{-1}=0$. The~element
$a_\epsilon$ is also irreducible and so $(a_\epsilon)_{-1}\gg0$. Now the~last equality yields
$a=b$, a~contradiction.

{\bf (b)} As was mentioned above, $A$ is a~Fredholm algebra and the~set $\Phi(A)$ of Fredholm
elements can be defined via the~formula~(\ref{27}). The~coefficients $a_{-1}$ and $a_0$ of
the~Laurent series expansion of $R(\cdot,a)$ around $r(a)$ satisfies ${\bf e}-(r(a)-a)a_0={\bf
e}-a_0(r(a)-a)=a_{-1}$ and, hence, $r(a)-a\in\Phi(A)$. Now the~required assertion follows at once
from part~{\bf (a)}.

{\bf (c)} It suffices to observe the~inclusion $a_{-1}\in{\cal F}(A)$ and to use part~~{\bf
(b)}.\hfill$\Box$

\section{The~closedness of the~center}\label{sec6}

A~normed algebra with a~unit ${\bf e}$ and with a~(closed, convex) cone $A^+$ is called an~{\it
ordered normed algebra}~\cite{Al3} if ${\bf e}\ge0$ and the~inequalities $a,b\ge0$ imply $ab\ge0$.
The~{\it center}~\cite{Al3} of an ordered normed algebra $A$ is called the~order ideal $A_{\bf e}$
generated by ${\bf e}$ (see Example~\ref{exm2}{\bf (b)}), i.e.,
$$A_{\bf e}=\{a\in A:-\lambda{\bf e}\le a\le\lambda{\bf e} \ \text{for some} \
\lambda\in{\Bbb R}^+\}.$$ The~aim of this section is to prove the~closedness of the~center $A_{\bf
e}$ in $A$ (Theorem~\ref{thm38}).

As was shown in \cite{Al3}, if an~element $a\in A_{\bf e}$ then $a^2\ge0$. The~next result makes
more precisely this fact.

\begin{lem}\
In an~arbitrary ordered normed algebra $A$ the~following identity holds
$$A_{\bf e}=\{a\in A:(\lambda+a)^2\ge0 \ \text{for every} \ \lambda\in{\Bbb R}\}.$$
\end{lem}

{\bf Proof.} In view of the~remarks above, it suffices to show that if $(\lambda+a)^2\in A^+$ for
all real $\lambda$ then $a\in A_{\bf e}$. To this end, we consider the~Banach algebra $B$ being
a~completion of $A$ and $K=\overline{A^+}$, where the~closure was taken in $B$. Obviously, $K$ is
a~wedge and $K\cdot K\subseteq K$. The~relations $a^2\in A^+$ and $(\lambda^2-a^2)^2\in A^+$ with
$\lambda\in{\Bbb R}$ are valid. For $\lambda>r(a)$ the~element $(\lambda^2-a^2)^{-1}$ is well
defined and belongs to $K$. Whence $\lambda^2-a^2\in A^+\cdot K\subseteq K$. On the~other hand,
$\lambda^2-a^2\in A$. In view of the~closedness of~$A^+$ in $A$, $\lambda^2-a^2\in K\cap A=A^+$.
Therefore, $0\le_{A^+}a^2\le_{A^+}\le\lambda^2{\bf e}$ and, in~particular, $a^2\in A_{\bf e}$.
Next, using the~relation $({\bf e}\pm a)^2\in A^+$, we obtain $\pm2a\le_{A^+}{\bf e}+a^2$. Finally,
$a\in A_{\bf e}$.\hfill$\Box$

\begin{thm}\label{thm38}\
The~center $A_{\bf e}$ is a~closed subset of an~arbitrary ordered normed algebra~$A$.
\end{thm}

{\bf Proof.} Consider a~sequence $\{a_n\}$ in the~center $A_{\bf e}$ satisfying $a_n\to a$ in $A$
as $n\to\infty$. Then $(\lambda+a_n)^2\ge0$ for all $\lambda\in{\Bbb R}$ and, hence,
$(\lambda+a)^2\ge0$. In view of the~preceding lemma, the~element $a\in A_{\bf e}$.\hfill$\Box$
\smallskip

The~Minkowski norm $\|\cdot\|_{\bf e}$ can be defined on the~center $A_{\bf e}$ (see~(\ref{28})).
Under this norm, $A_{\bf e}$ is a~(real) ordered normed algebra. The~next result follows
immediately from the~equality $A_{\bf e}=\bigcup\limits_{n=1}^\infty n[-{\bf e},{\bf e}]$,
the~Baire category theorem, and the~preceding theorem.

\begin{cor}\
In an~ordered Banach algebra $A$ the~embedding $(A_{\bf e},\|\cdot\|_A)\to(A_{\bf e},\|\cdot\|_{\bf
e})$ is continuous. In~particular, there exists a~constant $c>0$ satisfying $\|a\|_{\bf e}\le
c\|a\|_A$ for all $a\in A_{\bf e}$.
\end{cor}

\end{document}